\documentclass{article}

\usepackage{import}
\usepackage{macros}

\title{Density and holes of the visible points in certain cut-and-project sets}
\author{Hayim Shamir}
\date{}

\begin{document}
\maketitle

\begin{abstract}
We study the set of visible points from the origin in certain cut-and-project sets, specifically the density of the set. We generalize to all real quadratic number fields the formulas for the density found by G. Hammarhjelm, and the corresponding error bounds proved by I. Gringlaz, R. Kumar, B. Weiss. We also find the possible densities of the set of visible points in the cut-and-project set, with a constant generic lattice, and a changing window with constant volume. We also study bounds for the distance of large holes in the set of visible points of certain cut-and-project sets, generalizing theorems originally proved by P. Erd\H{o}s for $\mathbb{Z}^2$.
\end{abstract}

\section{Introduction}
Let $m \ge 1, d \ge 2$ be integers, and set $n = m + d$. We consider $\bR^n$ as the decomposition $\bR^n = \bR^m \oplus \bR^d$. We call $\bR^m$ the internal space and $\bR^d$ the physical space, with the natural projections $\pi_{\rm int}: \bR^n \rightarrow \bR^m$, $\pi_{\rm phys}: \bR^n \rightarrow \bR^d$.

For a locally finite point set $\cP \subset \bR^n$, and a bounded subset $\cW \subset \bR^m$ we define:
\[
\Lambda(\cW, \cP) = \{ \pi_{\rm phys}(\bm{y}) : \bm{y} \in \cP, \pi_{\rm int}(\bm{y}) \in \cW\}\,.
\]
This is also a locally finite point set $\Lambda(\cW, \cP) \subset \bR^d$. A {\em{window}} is a bounded subset $\cW \subset \bR^m$, with nonempty interior, and is Jordan measurable (boundary of measure 0).
When $\cW \subset \bR^m$ is a window, and $\cL \subset \bR^n$ is a lattice (a subset of the form $g \bZ^n$ for $g \in GL_n(\bR)$), the set $\Lambda(\cW, \cL)$ is called a {\em{cut-and-project set}} (also referred to as a model set). Throughout the paper, $\cL$ denotes a lattice in $\bR^n$, and $\cW \subset \bR^m$ denotes a window.
The cut and project set is called {\em{irreducible}} if $\pi_{\rm int}(\cL)$ is dense in $\bR^m$, and $\pi_{\rm phys}|_{\cL}$ is injective. Note that irreducibility depends only on the lattice $\cL$, and not on the window $\cW$.

Cut and project sets and their properties have been extensively studied (see \cite{BaakeGrimm2013,Marklof_2014,Meyer1972}), as they let us construct special locally finite point sets, especially non-periodic ones like the Penrose point set (see \cite[Theorem 3.3]{hammarhjelm2019densityminimalgapvisible}, originally proved in \cite{DEBRUIJN198153}), and the Ammann-Beenker point set (see \cite[Example 7.8]{BaakeGrimm2013}, originally proved in \cite{BeenkerNonPeriodicTilings1982}).

Let $D \subset \bR^d$ be a bounded set with $\vol(D) > 0$ and boundary of zero measure. We call such a set an {\em{averaging set}}, and in \Cref{sec:non_h_density,sec:basic_rings_different_density_example,sec:holes} we will also assume it is convex. The usual averaging sets taken are unit balls with respect to some metric.
Let $\cP \subset \bR^d$ be a locally finite point set, for $T > 0$, denote the dilation $TD = \{T\bm{x}: \bm{x} \in D\}$, and we define the $\textit{asymptotic density of } \cP \textit{ with respect to } D$:
\[
\theta(\cP) \overset{\text{def}}{=} \lim_{T\rightarrow\infty}{\frac{\# (\cP \cap TD)}{\vol(TD)}}\,.
\]
The existence of this limit and its value  might depend on the averaging set $D$. If the limit is defined and is the same for all averaging sets $D$, we say $\cP$ {\em{has density}}. We omit from the notation the dependence on the averaging set $D$, as in our cases, we will see below that the density will not depend on it.

It is well known that irreducible cut-and-project sets $\Lambda(\cW,\cL)$ have density given by the formula:
\begin{equation} \label{eq:irreducible_cut_and_project_density}
\theta\left( \Lambda(\cW, \cL) \right) = \frac{\vol(\cW)}{\covol(\cL)}\,,
\end{equation}
see \cite{hof1998uniform, Marklof_2014, schlottmann1998cut, Moody_2002}.

For any locally finite point set $\cP \subseteq \bR^d$, denote $\cP_{\star} = \cP \setminus \{0\}$, and define the set of visible points from the origin:
\[
\cP_{vis} = \{ \bm{x} \in \cP_{\star}\  :\ \forall t \in (0,1), t \bm{x} \notin \cP_{\star} \}\,.
\]
For a cut and project set, we denote 
\[
\Lambda_{\star}(\cW, \cL) = \Lambda(\cW,\cL)_{\star},\  \Lambda_{vis}(\cW, \cL) = \Lambda(\cW,\cL)_{vis}\ .
\]
In this paper we are interested in some properties of the set of visible points in certain cases of cut-and-project sets, most notably the density.

Denote by $\mu$ the Haar-measure on $GL_n(\bR)$, and we can look at the lattice defined by $\cL_g = g \bZ^n$, where $g$ is chosen outside of a set of zero $\mu$ measure. It is well known that for $\mu \textit{-a.e. } g$, cut-and-project sets with $\cL_g$ are irreducible. Indeed, both the failure of $\pi_{\rm int}(\cL_g)$ being dense in $\bR^m$, and the failure of injectivity of $\pi_{\rm phys}|_{\cL_g}$, can be seen to occur only on a countable collection of sub-manifolds of $GL_n(\bR)$ of positive co-dimension, so they have $\mu$-measure 0.

It was shown by Marklof and Str{\"o}mbergsson \cite{marklof2014visibilitydirectionsquasicrystals} that for any irreducible cut-and-project set, $\Lambda_{vis}(\cW,\cL)$ always has positive density, which doesn't depend on the averaging set, although without a formula for the density. In specific cases a formula is known, such as in \cite{hammarhjelm2019densityminimalgapvisible}, where specific lattices are considered, and it is also assumed that $\cW$ is star-shaped and for some constant $c > 1$, depending on which lattice is considered, we have $-\cW \subseteq c\cW$. Another example is \cite[Theorem 1.2]{gringlaz2025errorboundsvisiblepoints}, where it is proved that for $\mu$-a.e. $g$, assuming $\cW$ is star-shaped with respect to the origin, the density is $\theta(\Lambda_{vis}(\cW,\cL)) = \frac{1}{\zeta(n)}\theta(\Lambda(\cW,\cL))$. An additional noteworthy example is \cite{kumar2026simultaneous}, where the authors consider the set of points in $\Lambda(\cW,\cL)$ simultaneously visible from multiple different points instead of only the origin, and find formulas for the density for some similar cases.

One of the more interesting constructions of cut and project sets, which also lets us construct the Penrose point set and Ammann-Beenker point set, is using lattices from algebraic number fields. Let $K$ be a real quadratic number field, i.e. $K = \bQ(\sqrt{\fkD})$ for $\fkD > 1$ square-free. Let $\cO_K$ denote the ring of integers, and $\sigma(a+b\sqrt{\fkD}) = a-b\sqrt{\fkD}$ the nontrivial Galois automorphism of $K$. We have the standard Minkowski embedding of $\cO_K$ into $\bR^2$, defined by $\cL_{K}= \{(x,\sigma(x)) : x \in \cO_K \}$. With this we can construct the lattice:
\begin{equation} \label{eq:multi_dimension_number_field_lattice}
    \cL = \{(x_1, ..., x_d,\sigma(x_1), ..., \sigma(x_d)) :\  x_1,...,x_d \in \cO_K \} \cong \bigoplus^d\cL_{K} \subset \bR^{2d}\ .
\end{equation}
Here for our cut-and-project set, $m = d, n = 2d$. Also denote $\zeta_K$ the Dedekind zeta function of $K$, see \Cref{sec:non_h_density} for the definition.

Assuming the window $\cW \subseteq \bR^d$ is convex and symmetric around the origin, the density for the visible points in $\Lambda(\cW,\cL)$ was computed and proved by Hammarhjelm \cite{hammarhjelm2019densityminimalgapvisible}, and given error bounds by Gringlaz, Kumar and Weiss \cite{gringlaz2025errorboundsvisiblepoints}, but the proof only works for number fields $K$ satisfying the H-condition (see \Cref{def:Hammarhjelm_condition}). Rudnick \cite{rudnick2026classification} proved that only finitely many such number fields exist. Namely:
\begin{equation} \label{eq:full_H_condition_list}
\bQ(\sqrt{2}), \bQ(\sqrt{5}), \bQ(\sqrt{13}), \bQ(\sqrt{29}), \bQ(\sqrt{53}), \bQ(\sqrt{173}), \bQ(\sqrt{293}).
\end{equation}
In \Cref{sec:non_h_density} we find the density for all real quadratic number fields $K$, including convergence rates, thus generalizing the results from \cite{gringlaz2025errorboundsvisiblepoints,hammarhjelm2019densityminimalgapvisible}:
\begin{thm} \label{thm:main_thm}
    Given a real quadratic number field $K$, there exist constants $\alpha_1,...,\alpha_k \in \bQ, \beta_1,...,\beta_k \in K$, which depend only on $K$, for which $0 < \beta_i < \alpha_i \le 1$, such that for any $d \ge 2$, a convex and centrally symmetric window $\cW$ and a convex averaging set $D$ we have:
    \[
    \frac{\# (\Lambda_{vis}(\cW,\cL) \cap TD)}{\vol(TD)} = \sum_{i=1}^k{(\alpha_i^d - \beta_i^d)} \cdot \frac{\theta(\Lambda(\cW,\cL))}{\zeta_K(d)} + \begin{cases} 
    O\left(\frac{\log(T)}{\sqrt{T}}\right) & \text{if } d = 2, \\
    O\left(\frac{1}{\sqrt{T}}\right) & \text{if } d \ge 3.
    \end{cases}
    \]
    Where $\cL$ is defined by \eqref{eq:multi_dimension_number_field_lattice}, and the constants in the error term can depend on $\cW,D$, but not $T$.
\end{thm}

\begin{rmk}
    $\alpha_1,...,\alpha_k,\beta_1,...,\beta_k$ are not mysterious constants, we will provide an algorithm to find them for every real quadratic number field $K$. See \Cref{rmk:finding_the_constants} for the algorithm, and \Cref{rmk:constants_table} for the constants in the first few quadratic fields.
\end{rmk}

\begin{rmk} \label{rmk:H_condition_constants}
    For the quadratic fields $K$ in which H-condition is satisfied (see \Cref{def:Hammarhjelm_condition}), we have $k=1$ and the constants are just $\alpha_1 = 1, \beta_1 =\frac{1}{\lambda_K}$, where $\lambda_K > 1$ is the fundamental unit of $K$, see \Cref{prop:h_cond_visible_description}.
\end{rmk}

Next, we would like to discuss how the density $\theta(\Lambda_{vis}(\cW,\cL))$ can depend on the window $\cW$, even with the same lattice. In cases such as in the lattices arising from real quadratic number fields, for \Cref{thm:main_thm} we assume convexity of the window and central symmetry (in \cite{hammarhjelm2019densityminimalgapvisible} it was required that $-\cW \subseteq c \cW$ for $c > 1$ depending on $K$ instead of central symmetry). Another example is \cite[Theorem 1.2]{gringlaz2025errorboundsvisiblepoints}, with the lattice $\cL_g$ for $\mu$-a.e. $g$, where it is assumed that $\cW$ is star-shaped w.r.t. the origin.
In those formulas, the density depends only on the volume of $\cW$.

We show that in general when removing the convexity/star-shaped assumptions on $\cW$, the density does not depend solely on the volume of the window. In \Cref{sec:basic_rings_different_density_example} we give examples for two different windows $\cW_1, \cW_2 \subset \bR^m$ with the same volume, but $\theta(\Lambda_{vis}(\cW_1,\cL)) \ne \theta(\Lambda_{vis}(\cW_2,\cL))$, using the lattice from \eqref{eq:multi_dimension_number_field_lattice} with $K = \bQ(\sqrt2)$.

Afterwards, we ask a more general question: What is the possible range of the density $\theta(\Lambda_{vis}(\cW,\cL))$, when keeping the lattice and the volume of the window the same, but changing the window? In \Cref{sec:density_range_random_lattice}, we prove the following:
\begin{thm} \label{thm:random_density_range}
    For $\mu \textit{-a.e. } g \in GL_n(\bR)$, taking $\cL_g = g \bZ^n$, the following holds:

    \begin{enumerate}[label=(\arabic*)]
        \item For any window $\cW \subset \bR^m$, we have that:
        \[
        \frac{1}{\zeta(n)} \theta\left( \Lambda(\cW,\cL_g) \right) \le \theta\left( \Lambda_{vis}(\cW,\cL_g) \right) < \theta\left( \Lambda(\cW,\cL_g) \right).
        \]
        \item For every $V > 0$, and every $q \in \left[\frac{1}{\zeta(n)}\frac{V}{\covol(\cL_g)}, \frac{V}{\covol(\cL_g)}\right)$, there exists a window $\cW \subset \bR^m$ with $\vol(\cW) = V$, such that $\theta(\Lambda_{vis}(\cW,\cL_g)) =  q$.
    \end{enumerate}
\end{thm}
In fact, we prove this via a formula for the density of the visible points for $\mu \textit{-a.e. } g \in GL_n(\bR)$:
\begin{thm} \label{thm:random_density_formula}
    For $\mu \textit{-a.e. } g \in GL_n(\bR)$, taking $\cL_g = g \bZ^n$, for any window $\cW \subset \bR^m$, we have that:
        \[
        \theta\left(\Lambda_{vis}(\cW, \cL_g)\right) = \frac{1}{\zeta(n) \covol(\cL_g)}\sum_{k=1}^{\infty}{\frac{1}{k^n} \vol\left( \cW \setminus \left( \bigcup_{i=1}^{k-1}\frac{k}{i}\cW \right)\right)}.
        \]
\end{thm}
\begin{rmk}
    This gives us a formula for the density of the visible points, the sum is finite in cases such as when $\bm{0} \in \text{int}(\cW)$ or $\cW$ is star-shaped with respect to the origin, in the latter we are only left with the $k=1$ term, and this case and calculation of the density was already in \cite[Theorem 1.2]{gringlaz2025errorboundsvisiblepoints}.
\end{rmk}
\begin{rmk}
    In our proof of those theorems, we will get something stronger, we will show both \Cref{thm:random_density_range,thm:random_density_formula} hold for all $g \in GL_n(\bR)$ satisfying concrete conditions, mainly that $\cL_g$ gives irreducible cut-and-project sets, and another important condition, see \Cref{def:full_measure_generic_set}.
\end{rmk}

Finally, in \Cref{sec:holes}, we discuss a different property of the set of visible points- finding large holes with no points from the set $\Lambda_{vis}(\cW,\cL)$, specifically with lattices coming from number fields $K$ satisfying H-condition. Denote $D_0 = [0,1]^d$ for this part. For example, it was proved by Erd\H{o}s \cite{Erdos_1958} that $\bZ_{vis}^d$ has arbitrarily large holes, i.e. for every $T > 0$, there exists $\bm{x} \in \bR^d$, such that $(\bm{x} +TD_0) \cap \bZ_{vis}^d = \emptyset$. It was proven in \cite[Theorem 1.1]{gringlaz2025errorboundsvisiblepoints} that assuming $\cW$ is star-shaped with respect to the origin, there are holes as large as we want with no points in $\Lambda_{vis}(\cW,\cL)$, i.e. for any $T > 0$ exists $\bm{x} \in \bR^d$ such that $(\bm{x} +TD_0) \bigcap \Lambda_{vis}(\cW,\cL) = \es$. As was proved in \cite{Erdos_1958} for $\bZ_{vis}^n$, we prove bounds on how far those holes are from the origin, depending on the size:
\begin{thm} \label{thm:large_holes_distance}
    Let $K$ be a number field satisfying the H-condition, and $\cL \subset \bR^{2d}$ from \eqref{eq:multi_dimension_number_field_lattice}, and $\cW \subset\bR^d$ a convex and centrally symmetric window. 
     \begin{enumerate}[label=(\arabic*)]
        \item For all $\eps > 0$, for $T > 0$ large enough, we can find $\bm{x} \in \bR^d$ such that:
        \[
        \begin{aligned}
            &(\bm{x} +TD_0) \bigcap \Lambda_{vis}(\cW,\cL) = \es, \\
            &\|\bm{x}\| \le \exp{\left( (1 + \eps) \cdot \ \frac{\left(1 - \frac{1}{\lambda_K^d} \right) \vol(\cW)}{\covol(\cL)} d \cdot T^d \log(T)\right)}.
        \end{aligned}
        \]
        \item There exists $c > 0$, such that for all $T > 0$ large enough, if $\bm{x} \in \bR^d$ satisfies:
        \[
        (\bm{x} +TD_0) \bigcap \Lambda_{vis}(\cW,\cL) = \es,
        \]
        then we must have that:
        \[
        \|\bm{x}\| \ge \exp{\left( c T^{d-1} \log(T)\right)}.
        \]
    \end{enumerate}
    The logarithm and exponent are with respect to the natural base, $e$, and $\lambda_K > 1$ is the fundamental unit of the number field, see \Cref{sec:non_h_density}.
\end{thm}

\begin{rmk}
    If we normalize the cut and project set from the origin such that it would have density $1$, i.e. looking at $\Lambda' = \sqrt[d]{\frac{\vol(\cW)}{\covol(\cL)}} \Lambda(\cW,\cL)$, finding a hole $TD_0$ in $\Lambda_{vis}'$, is equivalent to finding a hole $\sqrt[d]{\frac{\covol(\cL)}{\vol(\cW)}} TD_0$ in $\Lambda_{vis}(\cW,\cL)$. Thus there exists an $\bm{x}$ with $(\bm{x} +TD_0) \bigcap \Lambda'_{vis} = \emptyset$ and:
    \[
        \|\bm{x}\| \le \exp{\left( (1 + \eps) \cdot \ \left(1 - \frac{1}{\lambda^d} \right) \cdot d T^d \log(T)\right)}.
    \]
    We can compare \Cref{thm:large_holes_distance}(1) to a past result from \cite{Erdos_1958}, where Erd\H{o}s gave an upper bound for the distance of holes in $\bZ_{vis}^2$ by using the Chinese remainder theorem, which can easily be generalized to $\bZ_{vis}^d$. The bound says that for any $\eps > 0$, for large enough $T > 0$, we can find $\bm{x} \in \bR^d$, such that $(\bm{x} +TD_0) \bigcap \bZ_{vis}^d = \emptyset$ and:
    \begin{equation} \label{eq:erdos_upper_bound}
        \|\bm{x}\| \le \exp{\left( (1 + \eps) \cdot d T^d \log(T)\right)}.
    \end{equation}
    We can see that in $\Lambda'$ we got a slightly lower bound, which comes from the additional blocking we can easily control in the proof, caused by multiplying by units of the number field.

    In \Cref{thm:large_holes_distance}(2), by tightening the bounds in the proof, one can take $c = \frac{\theta(\Lambda_{vis}(\cW,\cL))}{A \covol(\cL_K)} - \eps$ for any $\eps > 0$, where $A >0$ satisfies $\cW \subseteq [-A,A]^d$. Normalizing again, we get that if $\bm{x} \in \bR^d$ satisfies $(\bm{x} + TD_0) \bigcap \Lambda_{vis}' = \emptyset$, then:
    \begin{equation} \label{eq:lower_bound_holes_normalized}
    \| \bm{x}\| \ge \exp{\left((1-\eps)\frac{\sqrt[d]{\theta(\Lambda(\cW,\cL))}}{A\covol(\cL_K)} T^{d-1} \log(T)\right)}.
    \end{equation}
    We believe this bound is mainly due to our methods, and isn't optimal. It should be noted the upper bound and lower bound in the theorem are quite far apart, their exponents differ by a factor of $T$. Even for $\bZ_{vis}^2$, it is an open problem to close the gap between the bounds. The lower bound proved by Erd\H{o}s \cite{Erdos_1958}, was that if $\bm{x} \in \bR^2$ satisfies $(\bm{x} +TD_0) \bigcap \bZ_{vis}^2 = \emptyset$, then:
    \[
    \|\bm{x}\| \ge \exp{\left(\frac{1}{4}T \log{(T)} \right)}.
    \]
    So for $\bZ_{vis}^2$ as well, the lower bound differs from the upper bound \eqref{eq:erdos_upper_bound} (in the case $d=2$) by a factor of $T$ in the exponent.
\end{rmk}

\section{Visible density in cut and project sets for which H condition doesn't hold}\label{sec:non_h_density}
\subsection{Notations}

In this section, and in \Cref{sec:basic_rings_different_density_example,sec:holes} we will use standard notations and facts from algebraic number theory, see \cite{marcus2018number} for a standard reference.

Let $K$ be a real quadratic number field, i.e. there exists $\fkD>1$ which is square free such that $K=\bQ(\sqrt{\fkD})$.
Denote $\cO_K$ to be the ring of integers, $U_K = \cO_K^\times$ the unit group, and $C_K$ the class group. As the only roots of unity in $\cO_K$ are $\pm1$, by Dirichlet's unit theorem $U_K/\{\pm1\} \cong \bZ$, so we can find a fundamental unit $\lambda=\lambda_K>1$ such that $U_K=\{-1,1\} \times \{ \lambda^i : i \in\bZ \}$.
Denote by $\sigma(a+b\sqrt{\fkD}) = a-b\sqrt{\fkD}$ the non trivial Galois automorphism of $K$.

For an ideal $I \vartriangleleft \cO_K$, denote $N(I) = |\cO/I|$ to be the norm of the ideal in $\cO_K$, which is well known to be finite and multiplicative. Also for $g \in \cO_K$, denote $N(g) = N((g))$. We also have that $N(g) = |N_{K/\bQ}(g)| = |g\sigma(g)|$. And for $g = \frac{a}{b} \in K$ with $a,b \in \cO_K$, we also define $N(g) = \frac{N(a)}{N(b)} = |N_{K/\bQ}(g)| = |g\sigma(g)|$. Note for any $I \vartriangleleft \cO_K$, we have $N(I) \in \bZ_{\ge0}$, and it is well known there are only finitely many ideals in $\cO_K$ of each norm, as can be seen in the proof of \cite[Theorem 35, Corollary 2]{marcus2018number}, or by the more extensive \cite[Theorem 39]{marcus2018number} .

For $s \in \bR$ with $s > 1$, the Dedekind zeta function of the number field $K$ is given by:
\[
    \zeta_K(s) = \sum_{0 \ne I \vartriangleleft \cO_K}{\frac{1}{N(I)^s}}\ .
\]
If we take the M\"obius function given by:
\[
    \mu_K(I) = \begin{cases} 
    1 & \text{if } I = (1), \\
    (-1)^r & \text{if } I = P_1 \cdots P_r \text{ for distinct prime ideals } P_1,...,P_r, \\
    0 & \text{if } I\subset P^2 \text{ for some prime ideal } P,
    \end{cases}
\]
we have the identity:
\[
    \frac{1}{\zeta_K(s)} = \sum_{0 \ne I \vartriangleleft \cO_K}{\frac{\mu_K(I)}{N(I)^s}}\ .
\]

Denote $\cL_{K}= \{(x,\sigma(x)) : x \in \cO_K \}$ as the standard Minkowski embedding of $\cO_K$. We denote the discriminant $D_K=\covol(\cL_{K})^2$ (it is positive because $K$ is real and quadratic).

A standard fact proved by Minkowski, is that for any ideal class in $C_K$, there exists an ideal $0 \ne J \vartriangleleft \cO_K$ in that ideal class, with norm bounded by:
\begin{equation} \label{eq:minkowski_class_representive_bound}
    N(J) \le \frac{n!}{n^n} \left(\frac{4}{\pi}\right)^s \sqrt{|D_K|}\ .
\end{equation}
Where $n = [K:\bQ]$ is the extension degree, and $s$ is the number of pairs of complex embeddings (For field embeddings $\varphi: K \rightarrow \bC$ such that $\Img{\varphi} \nsubseteq \bR$, we say $\{ \varphi, \bar{\varphi} \}$ is a pair of complex embeddings). In our case, $n = 2, s = 0, D_K > 0$, and the resulting bound is $N(J) \le \frac{1}{2} \sqrt{D_K}$. See \cite[Theorem 37, Corollary 2]{marcus2018number}.

For $d \ge 2$, we define the lattice $\cL$ as in \eqref{eq:multi_dimension_number_field_lattice}.
We also extend $\sigma$ to $K^d$ by defining $\sigma(x_1,...,x_d)=(\sigma(x_1),...,\sigma(x_d))$.

Given a window $\cW \subset \bR^{d}$, we define the cut-and-project set $\Lambda(\cW,\cL)$, and the set of visible points $\Lambda_{vis}(\cW,\cL)$ as before. We use the following terminology throughout the proof: we say a point $0 \ne \bm{x} \in \bR^d$ {\emph{blocks}} the point $0 \ne \bm{y} \in \bR^d$ if there exists $t \in (0,1)$ such that $\bm{x} = t \bm{y}$. The set of visible points $\Lambda_{vis}(\cW,\cL)$ is the set of points in $\Lambda_{\star}(\cW,\cL)$ not blocked by any other point in $\Lambda_{\star}(\cW,\cL)$.

For an ideal $0 \ne I \vartriangleleft \cO_K$, we define the sub-lattice $\cL_I \le \cL$:
\[
    \cL_I = \{ (x_1,...,x_d,\sigma(x_1),...,\sigma(x_d)) \in \cL : x_1,....,x_d \in I\}\,.
\]
We have that $\covol(\cL_I) = N(I)^d \covol(\cL)$.

For $0 \ne \fkg \in K$, let $a_\fkg$ be the diagonal matrix:
\[
a_\fkg \overset{\text{def}}{=} \text{diag}( \underbrace{\fkg,...,\fkg}_{\textit{d times}}, \underbrace{\sigma(\fkg),...,\sigma(\fkg)}_{\textit{d times}} ).
\]
Note that for an ideal $0 \ne I \vartriangleleft \cO_K$, if $\fkg I \vartriangleleft \cO_K$, we have that $\cL_{\fkg I} = a_\fkg \cL_I$, and in particular if $\fkg \in U_K$ is a unit, every $\cL_I$ is invariant under $a_\fkg$. Also note for ideals $0 \ne I,J \vartriangleleft \cO_K$, there exists $0 \ne \fkg \in K$ such that $a_\fkg \cL_I = \cL_J$ if and only if $I,J$ are of the same class in the class group.

And for a window $\cW \subset \bR^{d}$ and ideal $0 \ne I \vartriangleleft \cO_K $, we define the $I$-primitive points:
\[
\Lambda_{prim}(\cW,\cL,I) = \{ (x_1,...,x_d) \in \Lambda(\cW, \cL_I) : \gcd(x_1,...,x_d) = I \}\,.
\]
For an ideal $0 \ne I \vartriangleleft \cO_K$ and a point $\bm{x} = (x_1,...,x_d) \in \cO_K^d$, we say $I \textit{ divides } \bm{x}$ if $I\ |\ \gcd(x_1,...,x_d)$, or in other words, $x_1,...,x_d \in I$ (using the unique ideal factorization thanks to $\cO_K$ being a Dedekind domain).

Throughout this section, and \Cref{sec:holes} we will have additional assumptions on the window and the averaging set:
\begin{defn}
    A window $\cW \subset \bR^m$ is called a nice window, if in addition to the previous requirements of a window, it is also convex and centrally symmetric.
    An averaging set $D \subset \bR^d$ is called a nice averaging set, if in addition to the previous requirements of an averaging set, it is also convex.
\end{defn}
Note every nice window is also star-shaped with respect to the origin.
\begin{rmk}
    The assumption of convexity from the averaging set is not usually made in the literature, we need it mainly for the error rates, for \Cref{lem:point_counting_orbit_L_I}. Many other averaging sets can be approximated using convex sets, or taking disjoint unions / subtraction of a convex subset.
    
    In addition, by \cite{marklof2014visibilitydirectionsquasicrystals}, the visible points have density, and it's the same for any averaging set. So after we found the density for nice averaging sets here, the density of $\Lambda_{vis}(\cW,\cL)$ along any other averaging set must be the same, but maybe without the error rate we found here.
\end{rmk}

\begin{defn} \label{def:Hammarhjelm_condition}
    We say a real quadratic number field $K$ has the {\emph{Hammarhjelm condition}} (or H-condition for short), if $\cO_K$ is a principal ideal domain, and:
    \[
    \cL_K \cap ((1,\lambda_K) \times [-1,1] ) = \emptyset.
    \]
\end{defn}
\begin{rmk}
    As said earlier, this condition was used in both \cite{hammarhjelm2019densityminimalgapvisible,gringlaz2025errorboundsvisiblepoints} in order to calculate the density of $\Lambda_{vis}(\cW,\cL)$, assuming $K$ satisfies this property. The proof goes through a simple description for $\Lambda_{vis}(\cW,\cL)$ in those cases, see \Cref{prop:h_cond_visible_description}.
    However the proof from those papers doesn't work for other real quadratic fields, and Rudnick \cite{rudnick2026classification} proved there are finitely many number fields satisfying H-condition, see \eqref{eq:full_H_condition_list}.
\end{rmk}
\begin{prop}[{see \cite[Lemma 5.6]{gringlaz2025errorboundsvisiblepoints}}]
\label{prop:h_cond_visible_description}
    If $K$ satisfies the H-condition, for a nice window $\cW$ and the lattice $\cL \subset \bR^{2d}$ from \eqref{eq:multi_dimension_number_field_lattice}, we have:
    \[
    \Lambda_{vis}(\cW,\cL) = \Lambda_{prim}\left(\cW , \cL, (1) \right) \setminus \Lambda_{prim}\left(\frac{1}{\lambda}\cW , \cL, (1) \right).
    \]
\end{prop}

From here throughout the section, we will work with a constant quadratic number field $K$, constant nice averaging set $D$, and a constant nice window $\cW$, that whatever we call constant is allowed to depend on.

\subsection{Generalizing lemmas from \cite{gringlaz2025errorboundsvisiblepoints}}

\begin{lem} \label{lem:point_counting_orbit_L_I}
    Let $D \subset \bR^d$ be a nice averaging set, $\cW$ a nice window, $r > 0$ and $0 \ne J \vartriangleleft \cO_K$. Then there exists a constant $C = C(D, \cW, r, J)$ such that for any $0 \ne \fkg \in K$ and any $T \ge N(\fkg)$ we have:
    \[
    \left| \# (\Lambda(r \cW,a_\fkg \cL_J) \cap TD) - \frac{\vol(TD) \cdot \vol(\cW)}{\covol(\cL_J)} \cdot \frac{r^d}{N(\fkg)^d} \right| \le C \left( \frac{T}{N(\fkg)} \right)^{d-\frac{1}{2}}\,.
    \]
\end{lem}
\begin{proof}
    This is a generalization of \cite[Lemma 5.4]{gringlaz2025errorboundsvisiblepoints}, it is stated here for general $\fkg \in K$ instead of $\fkg \in \cO_K$, and $\cL_J$ instead of $\cL$, and in \cite{gringlaz2025errorboundsvisiblepoints} it is also assumed that $K$ satisfies H-condition. The proof from \cite{gringlaz2025errorboundsvisiblepoints} for the special case gives this more general statement, as $\cL_J$ is also invariant under $a_{\fkg_0}$ for $\fkg_0 \in U_K$, it is not used that $\fkg \in \cO_K$, and also H-condition is not used in the proof.
\end{proof}

\begin{lem} \label{lem:point_counting_all_L_I}
   Let $D \subset \bR^d$ be a nice averaging set, $\cW$ a nice window and $r > 0$. Then there exists a constant $C = C(D, \cW, r)$ such that for any $0 \ne I \vartriangleleft \cO_K$ and any $T \ge N(I)$ we have:
    \[
    \left| \# (\Lambda(r \cW,\cL_I) \cap TD) - \frac{\vol(TD) \cdot \vol(\cW)}{\covol(\cL)} \cdot \frac{r^d}{N(I)^d} \right| \le C \left( \frac{T}{N(I)} \right)^{d-\frac{1}{2}}\,.
    \]
\end{lem}
\begin{rmk}
    The main difference between \Cref{lem:point_counting_orbit_L_I} and this lemma is that here the constant doesn't depend on $I$, and $I$ here actually plays the role of $\fkg J$ from \Cref{lem:point_counting_orbit_L_I}. We get a uniform error bound for counting points in $\# (\Lambda(r \cW,\cL_I) \cap TD)$ over all the ideals $I$, in contrast to \Cref{lem:point_counting_orbit_L_I} which only got a uniform bound over all $\fkg J$ for a constant ideal $J$ and changing $0 \ne \fkg \in K$. Those are different, as in our setting, $\cO_K$ is not necessarily a principal ideal domain.
\end{rmk}
\begin{proof}
    Take representatives of each class in the class group $0 \ne J_1,...,J_h \vartriangleleft \cO_K$. Apply \Cref{lem:point_counting_orbit_L_I} to each of them, and denote the constants $C_1,...,C_r$. So we define $C = \max\limits_{i = 1,...,h}(C_i N(J_i)^{d-1/2})$. Now for any ideal $0 \ne I \vartriangleleft \cO_K$, we can take its corresponding class representative $J_i$, and $\frac{I}{J_i}$ must be a principal ideal, so we take a generator $(\fkg) = \frac{I}{J_i}$. So for any $T \ge N(I)$ we also have $T \ge N(\fkg)$, thus:
    \[
        \begin{aligned}
            &\left| \# (\Lambda(r \cW,\cL_I) \cap TD) - \frac{\vol(TD) \cdot \vol(\cW)}{\covol(\cL)} \cdot \frac{r^d}{N(I)^d} \right| = \\
            &= \left| \# (\Lambda(r \cW,a_\fkg \cL_{J_i}) \cap TD) - \frac{\vol(TD) \cdot \vol(\cW)}{\covol(\cL_{J_i})} \cdot \frac{r^d}{N(\fkg)^d} \right| \le \\
            &\le C_i \left( \frac{T}{N(\fkg)} \right)^{d-\frac{1}{2}} \overset{\star}{\le} C \left( \frac{T}{N(I)} \right)^{d-\frac{1}{2}}.
        \end{aligned}
    \]
    Where $\star$ follows from the choice of $C$. Thus we get as we wanted.
\end{proof}

\begin{lem} \label{lem:point_bounding_orbit_L_I}
    Let $D \subset \bR^d$ be a nice averaging set, $\cW$ a nice window, $r > 0$ and $0 \ne J \vartriangleleft \cO_K$. Then there exists a constant $M = M(D, \cW, r, J)$ such that for any $0 \ne \fkg \in K$ and any $T > 0$ we have:
    \[
    \# (\Lambda(r \cW,a_\fkg \cL_J) \cap TD) \le \frac{M \cdot T^d}{N(\fkg)^d}\,.
    \]
\end{lem}
\begin{proof}
    This is a generalization of \cite[Lemma 5.5]{gringlaz2025errorboundsvisiblepoints}, where the authors reference \cite[Lemma 4.6]{hammarhjelm2019densityminimalgapvisible}. As before the proof stays the same, it doesn't use that $\fkg \in \cO_K$, $\cL_J$ is also invariant under $a_{\fkg_0}$ for any $\fkg_0 \in U_K$, and it is not used that $K$ satisfies H-condition.
\end{proof}

\begin{lem} \label{lem:point_bounding_all_L_I}
    Let $D \subset \bR^d$ be a nice averaging set, $\cW$ a nice window. Then there exists a constant $L = L(r, D, \cW)$ such that for any $0 \ne I \vartriangleleft \cO_K$ and any $T > 0$ we have:
    \[
    \# (\Lambda(r \cW,\cL_I) \cap TD) \le \frac{L \cdot T^d}{N(I)^d}\,.
    \]
\end{lem}
\begin{proof}
    Take representatives for the class group $0 \ne J_1,...,J_h \vartriangleleft \cO_K$. Apply \Cref{lem:point_bounding_orbit_L_I} to each of them, and denote the constants $L_1,...,L_r$. So we define $L = \max\limits_{i = 1,...,h}(L_i N(J_i)^d)$. Now for any ideal $0 \ne I \vartriangleleft \cO_K$, we can take its corresponding class representative $J_i$, and $\frac{I}{J_i}$ must be a principal ideal, so we take a generator $(\fkg) = \frac{I}{J_i}$. So for any $T > 0$ we have:
    \[
        \# (\Lambda(r \cW,\cL_I) \cap TD) = \# (\Lambda(r \cW,a_\fkg \cL_{J_i}) \cap TD) \le \frac{L_i \cdot T^d}{N(\fkg)^d} \overset{\star}{\le} \frac{L \cdot T^d}{N(I)^d}\,,
    \]
    where $\star$ follows from the choice of $L$. Thus we get as we wanted.
\end{proof}

\subsection{Description of the visible points}

\begin{lem} \label{lem:finitely_bad_ideals}
    Given a nice window $\cW$, and $\cL \subset \bR^{2d}$ from \eqref{eq:multi_dimension_number_field_lattice}, there is a finite number of ideals $I$ which divide points in $\Lambda_{vis}(\cW,\cL)$.
    In particular, if $N(I) \ge \frac{\lambda}{2}\sqrt{D_K}$, then it doesn't divide any point $\bm{x} \in \Lambda_{vis}(\cW,\cL)$.
\end{lem}
\begin{proof}
    Because there are finitely many ideals of each norm, the second statement proves the first one.
    Let $0 \ne I \vartriangleleft \cO_K$ with $N(I) \ge \frac{\lambda}{2}\sqrt{D_K}$ and by way of contradiction, assume it divides $\bm{x} \in \Lambda_{vis}(\cW,\cL)$. $N(I) > 0$, so $I \ne (0)$.
    Take an ideal $0 \ne J \vartriangleleft \cO_K$ which is in the same ideal class as $I \vartriangleleft \cO_K$, with norm bounded by:
    \[
    N(J) \le \frac{2!}{2^2} \left( \frac{4}{\pi} \right)^0 \sqrt{D_K} = \frac{1}{2} \sqrt{D_K}\,.
    \]
    As in \eqref{eq:minkowski_class_representive_bound}. Thus $\frac{I}{J}$ must be a principal fractional ideal, so there exists $t \in K$ such that $\frac{I}{J} = \left( t \right)$. As $-1 \in \cO_K$, every ideal is equal to its minus, so we can assume without loss of generality that $t>0$.
    We will look at the point $\bm{y} = \lambda^\ell t^{-1} \cdot \bm{x}$, and show that there exists $\ell \in \bZ$, such that $\bm{y} \in \Lambda(\cW,\cL)$, and $\bm{y}$ blocks $\bm{x}$. 
    
    First, because $I$ divides $\bm{x}$, for all $1 \le i \le d$, the ideal $(y_i) = \frac{J}{I} (x_i)$ is integral, so we have $y_i \in \cO_K$. So $(\bm{y},\sigma(\bm{y})) \in \cL$.
    
    Second, in order to have $\bm{y} \in \Lambda(\cW,\cL)$, we need to have $\sigma(\bm{y}) \in \cW$. Note that $\sigma(\bm{y})=\sigma(\lambda^\ell t^{-1}) \sigma(\bm{x})$, thus because $\cW$ is centrally symmetric and star-shaped, and $\sigma(\bm{x}) \in \cW$, to get $\sigma(\bm{y}) \in \cW$, it is enough that $|\sigma(\lambda^\ell t^{-1})| \le 1$.
    
    And thirdly, we need $\bm{y}$ to block $\bm{x}$, so we need $0 < \lambda^\ell t^{-1} < 1$.
    
    So we get two requirements from $\ell$ which are:
    \[
    \begin{aligned} 
    &1. \ \ \  1 \ge |\sigma(\lambda^\ell t^{-1})| = \lambda^{-\ell} |\sigma(t)^{-1}|\,. \\
    &2. \ \ \  \lambda^\ell t^{-1} < 1\,.
    \end{aligned}
    \]
    This is equivalent to:
    \begin{equation} \label{eq:ell_bounds_finitely_many_bad_ideals}
        |\sigma(t)|^{-1} \le \lambda^\ell < t.
    \end{equation}
    And we have $\ell \in \bZ$ for which this holds, because the ratio between the two sides of \eqref{eq:ell_bounds_finitely_many_bad_ideals} satisfies:
    \[
    \frac{t}{|\sigma(t)|^{-1}} = |t\sigma(t)| = N(t) = \frac{N(I)}{N(J)} \overset{\star}{\ge} \lambda\,.
    \]
    Where $\star$ follows from the requirement of $N(I) \ge \frac{\lambda}{2}\sqrt{D_K}$, and the choice of $J$ with $N(J) \le \frac{1}{2} \sqrt{D_K}$. Thus for this $\ell$, $\bm{y} \in \Lambda(\cW,\cL)$, and it blocks $\bm{x} \in \Lambda_{vis}(\cW,\cL)$, a contradiction.
\end{proof}

So we have the finite set of ideals $I_1,...,I_k$ (including the ideal $(1)$) with norm at most $\frac{\lambda}{2}\sqrt{D_K}$, and any other ideal can't divide visible points in $\Lambda_{vis}(\cW,\cL)$. Thus we can write the visible points as a finite disjoint union:
\begin{equation} \label{eq:non_h_disjoint_union}
    \Lambda_{vis}(\cW,\cL) = \biguplus_{i=1}^k \left( \Lambda_{vis}(\cW,\cL) \cap \Lambda_{prim}(\cW,\cL,I_i) \right).
\end{equation}

\begin{lem} \label{lem:visible_points_classification}
    For every ideal $0 \ne I \vartriangleleft \cO_K$, there exists a constant $0 < r_I \le 1$ depending only on $K,I$, such that:
    \[
    \Lambda_{vis}(\cW,\cL) \cap \Lambda_{prim}(\cW,\cL,I) = \Lambda_{prim}(\cW,\cL,I) \setminus \Lambda_{prim}(r_I\cW,\cL,I).
    \]
\end{lem}

\begin{proof}
    We prove the contrapositive, that there exists $0 < r_I \le 1$ such that for $\bm{x} \in \Lambda_{prim}(\cW,\cL,I)$ we have:
    \begin{equation} \label{eq:invisibility_equivalence}
        \bm{x} \in \Lambda_{prim}(r_I\cW,\cL,I) \iff \bm{x} \notin \Lambda_{vis}(\cW,\cL).
    \end{equation}
    We will define $r_I$, and then show why this holds. Go over all of the divisors $J_1 | I$, and all possibilities for $J_2 \in \{I_1,...,I_k\}$ (the ideals which might divide visible points), such that $J_1,J_2$ are of the same class in the class group. For each such pair $(J_1,J_2)$, take an arbitrary $t \in K$ such that $(t) = \frac{J_2}{J_1}$, and also $t > 0$. Write $l = \lfloor \log_{\lambda}{t} \rfloor + 1$. We define:
    \begin{equation} \label{eq:maximum_for_r}
        r_I = \min\left(1, \max_{\text{all choices of } J_1,J_2}\left(\lambda^{-l} |\sigma(t)|^{-1}\right)\right).
    \end{equation}
    
    For the $\Leftarrow$ direction in \eqref{eq:invisibility_equivalence}, if $r_I = 1$ it is obvious. Otherwise, let $\bm{y} = s \bm{x}$ for $0 < s < 1$ be a point that blocks $\bm{x}$. As $\Lambda(\cW,\cL)$ is locally finite, we can also assume it is visible, $\bm{y} \in \Lambda_{vis}(\cW,\cL)$. Because of the first coordinate, we have $s \in K$, so there exists $a,b \in \cO_K$ such that $s = \frac{a}{b}$. Write $J_3 = \gcd(a,b)$, and $J_1 = \frac{(b)}{J_3}$, $J_2 = \frac{(a)}{J_3}$. For all $1 \le i \le d$, we have $J_1 (y_i) = J_2 (x_i)$. Also $\gcd(J_1,J_2) = 1$, thus we must have that $J_1$ divides each $x_i$, and $J_2$ divides each $y_i$, i.e. $J_1$ divides $\bm{x}$, and $J_2$ divides $\bm{y}$. Because $\bm{y} \in \Lambda_{vis}(\cW,\cL)$, we must have $J_2 \in \{ I_1,...,I_k \}$, and because $\bm{x} \in \Lambda_{prim}(\cW,\cL,I)$, we must have $J_1|I$.
    Also $J_1,J_2$ are in the same class, because $(s) = \frac{J_2}{J_1}$.
    
    Thus $J_1,J_2$ are one of the pairs we took the maximum over, let $t > 0$ denote the generator taken for $(t) = \frac{J_2}{J_1}$ in \eqref{eq:maximum_for_r}, thus with $l = \lfloor \log_{\lambda}{t} \rfloor + 1$, we have $\lambda^{-l} |\sigma(t)|^{-1} \le r_I$. Because $s,t \in K$ generate the same fractional ideal, $\frac{t}{s}$ must be a unit, and because they are both positive, we get there exists $q \in \bZ$ such that $t = \lambda^q s$. Because $t \lambda^{-q} = s < 1$, we get that $q > \log_{\lambda}{t}$, and because $q \in \bZ$, this means $q \ge \lfloor \log_{\lambda}{t} \rfloor + 1 = l$.
    Now we can see we have that:
    \[
    \sigma(\bm{x}) = \sigma(s^{-1} \bm{y}) = \sigma(s)^{-1} \sigma(\bm{y}) = \sigma(\lambda^{-q} t)^{-1} \sigma(\bm{y}) = \pm \lambda^{-q} \sigma(t)^{-1} \sigma(\bm{y}).
    \]
    And because $\bm{y} \in \Lambda(\cW,\cL)$, we must have $\sigma(\bm{y}) \in \cW$, thus because $\cW$ is star shaped and centrally symmetric:
    \[
    \sigma(\bm{x}) = \pm \lambda^{-q} \sigma(t)^{-1} \sigma(\bm{y}) \in (\lambda^{-q} |\sigma(t)|^{-1}) \cW \subseteq (\lambda^{-l} |\sigma(t)|^{-1}) \cW \subseteq r_I \cW.
    \]
    And $\bm{x}$ is primitive for $I$, so $\bm{x} \in \Lambda_{prim}(r_I\cW,\cL,I)$, and we finished the first direction.

    For the $\Rightarrow$ direction in \eqref{eq:invisibility_equivalence}, let $\bm{x} \in \Lambda_{prim}(r_I\cW,\cL,I)$, take ideals attaining the maximum $r_I$ in \eqref{eq:maximum_for_r}, $J_1,J_2$ and their corresponding $t, l$, as in the definition of $r_I$, so we have $r_I \le \lambda^{-l} |\sigma(t)|^{-1}$. We claim that $\bm{y} = \lambda^{-l} t \bm{x} \in \Lambda(\cW,\cL)$, and blocks $\bm{x}$. For each $1 \le i \le d$, we have $(y_i) = (t)(x_i) = J_2 \frac{(x_i)}{J_1}$, so because $J_1 | I | \bm{x}$, we have that $y_i \in \cO_K$, so $(\bm{y},\sigma(\bm{y})) \in \cL$.
    We also have $\sigma(\bm{y}) = \pm \lambda^l \sigma(t) \sigma(\bm{x})$, so because $\sigma(\bm{x}) \in r_I \cW$ and $\cW$ is star-shaped and centrally symmetric, we get:
    \[
    \sigma(\bm{y}) = \pm \lambda^l \sigma(t) \sigma(\bm{x}) \in \lambda^l |\sigma(t)| r_I \cW \subseteq \lambda^l |\sigma(t)| \cdot \lambda^{-l} |\sigma(t)|^{-1} \cW = \cW.
    \]
    thus $\bm{y} \in \Lambda(\cW,\cL)$. Finally, we need to explain why $\bm{y}$ blocks $\bm{x}$, i.e. why $0 < \lambda^{-l} t < 1$. It is obviously positive, and $t < \lambda^l$ because $l = \lfloor \log_{\lambda}{t} \rfloor + 1$. Thus $\bm{y} \in \Lambda(\cW,\cL)$ blocks $\bm{x}$, so $\bm{x} \notin \Lambda_{vis}(\cW,\cL)$ as we wanted.
\end{proof}

\begin{rmk} \label{rmk:finding_the_constants}
    Notice that given a quadratic number field $K$, we can label the ideals of norm at most $\frac{\lambda}{2} \sqrt{D_K}$ by $I_1,...,I_k$. For each $I_i$ we can apply \Cref{lem:visible_points_classification}, and find the corresponding $r_{I_i}$ by \eqref{eq:maximum_for_r} (we denote for short $r_i := r_{I_i}$). This will allow us to find the exact $\alpha_1,...,\alpha_k,\beta_1,...,\beta_k$ in the main theorem. They will be $\alpha_i = \frac{1}{N(I_i)}, \beta_i = \frac{r_i}{N(I_i)}$ (Note that some of $I_1,...,I_k$ might not divide visible points, in the cases when $r_i = 1$).
\end{rmk}

\begin{rmk}
    For \Cref{lem:finitely_bad_ideals,lem:visible_points_classification}, which give us the classification of the visible points, we didn't need convexity of the window $W$, but only the bounded, centrally symmetric and star-shaped with respect to the origin assumptions. The other assumptions are used in the point counting lemmas.
\end{rmk}

\subsection{Density calculation}

Let $I_1,...,I_k$ be the ideals with norm at most $\frac{\lambda}{2} \sqrt{D_K}$. For each $I_i$ we can apply \Cref{lem:visible_points_classification}, and find the corresponding $r_{I_i}$ by \eqref{eq:maximum_for_r} (we denote for short $r_i := r_{I_i}$). So by \eqref{eq:non_h_disjoint_union} and \Cref{lem:visible_points_classification}, we have:
\[
    \begin{aligned}
        &\Lambda_{vis}(\cW,\cL) = \biguplus_{i=1}^k \left( \Lambda_{vis}(\cW,\cL) \cap \Lambda_{prim}(\cW,\cL,I_i) \right) = \\
        &= \biguplus_{i=1}^k \left( \Lambda_{prim}(\cW,\cL,I_i) \setminus \Lambda_{prim}(r_i\cW,\cL,I_i) \right).
    \end{aligned}
\]
So because $\cW$ is star shaped:
\[
\#(\Lambda_{vis}(\cW,\cL) \cap TD) = \sum_{i=1}^k{\left( \#(\Lambda_{prim}(\cW,\cL,I_i) \cap TD) - \#(\Lambda_{prim}(r_i\cW,\cL,I_i) \cap TD) \right)}.
\]
So for the rest of the proof, it will be enough to count how many points are in $\Lambda_{prim}(r \cW,\cL,I)$ for $0 < r \le 1$, and we will sum it over a constant number of ideals and values of $r$, so it will not affect the error bounds.

\begin{lem} \label{lem:inclusion_exclusion}
    For any ideal $0 \ne I \vartriangleleft \cO_K, r > 0$, and $T > 0$ we have:
    \[
    \# (\Lambda_{prim}(r \cW, \cL, I) \cap TD) = \sum_{0 \ne J \vartriangleleft \cO_K}{\left( \mu_K(J) \cdot \# (\Lambda_{\star}(r \cW, \cL_{J \cdot I}) \cap TD \right)}.
    \]
\end{lem}
\begin{proof}
    This is as in \cite[Formula 24]{gringlaz2025errorboundsvisiblepoints}, follows easily from the local finiteness of $\Lambda(r \cW, \cL_I)$ using inclusion-exclusion.
\end{proof}

\begin{lem} \label{lem:primitive_point_counting}
    Let $0 \ne I \vartriangleleft \cO_K$ and $r > 0$. We have that:
    \[
    \begin{aligned}
        &\left| \# (\Lambda_{prim}(r \cW, \cL, I) \cap TD) - \frac{r^d}{N(I)^d} \cdot\frac{\vol(\cW) \cdot \vol(TD)}{\covol(\cL)} \cdot  \frac{1}{\zeta_K(d)} \right| \le \\
        &\le \begin{cases}
            C_2 T^{3/2} \log(T) & \text{if } d = 2,\\
            C_d T^{d-1/2} & \text{if } d > 2.
        \end{cases}
    \end{aligned}
    \]
    Here $C_2, C_d$ are constants which don't depend on $T$.
\end{lem}
\begin{proof}
    Exactly the same as in \cite[Lemma 5.7]{gringlaz2025errorboundsvisiblepoints}, with the new lemmas we have instead:
    We use \Cref{lem:inclusion_exclusion} instead of \cite[Equation 24]{gringlaz2025errorboundsvisiblepoints} for the inclusion-exclusion formula, \Cref{lem:point_counting_all_L_I} instead of \cite[Lemma 5.4]{gringlaz2025errorboundsvisiblepoints} to approximate the first terms in the series, which come from ideals $J \vartriangleleft \cO_K$ with $N(J\cdot I) \le T$, and finally \Cref{lem:point_bounding_all_L_I} instead of \cite[Lemma 5.5]{gringlaz2025errorboundsvisiblepoints} to bound the tail of the series, the terms which come from ideals $J \vartriangleleft \cO_K$ with $N(J \cdot I) > T$.
\end{proof}

\begin{proof} [Proof of \Cref{thm:main_thm}]
    Now summing the constants from \Cref{lem:primitive_point_counting}, when using it for each $I_1,...,I_k$ and both $r \in \{ 1, r_i\}$, we get in total that:
    \[
    \begin{aligned}
        &\frac{\# (\Lambda_{vis}(\cW,\cL) \cap TD)}{\vol(TD)} = \\
        &= \sum_{i=1}^k{\left( \frac{\#(\Lambda_{prim}(\cW,\cL,I_i) \cap TD)}{\vol(TD)} - \frac{\#(\Lambda_{prim}(r_i\cW,\cL,I_i) \cap TD)}{\vol(TD)} \right)} = \\
        &= \sum_{i=1}^k{\left( \frac{1}{N(I_i)^d} - \frac{r_i^d}{N(I_i)^d}\right)} \cdot \frac{\theta(\Lambda(\cW,\cL))}{\zeta_K(d)} + \begin{cases} 
        O\left(\frac{\log(T)}{\sqrt{T}}\right) & \text{if } d = 2, \\
        O\left(\frac{1}{\sqrt{T}}\right) & \text{if } d \ge 3.
        \end{cases}
    \end{aligned}
    \]
    So taking $\alpha_i = \frac{1}{N(I_i)}, \beta_i = \frac{r_i}{N(I_i)}$, (and removing the ones for which $r_i = 1$ or equivalently $\alpha_i = \beta_i$) we get the main theorem.
\end{proof}

\begin{rmk} \label{rmk:constants_table}
    As we said, the $\alpha_i, \beta_i$  from the main theorem, only depend on the quadratic field $K$. The following is a table of the corresponding ideals and constants for the first few quadratic fields. We omit ideals we found at the start when bounding the norm in \Cref{lem:finitely_bad_ideals} for which $r_i = 1$ and $\alpha_i=\beta_i$.
    \begin{table}[H]
    \centering
    \caption{Values of the constants for the first few Quadratic fields}
    \begin{tabular}{p{1cm} p{1cm} p{3cm} p{2cm} p{3cm}}
        \toprule
        \textbf{$K$} & \textbf{Has H property?} & \textbf{Ideals $I_i$ dividing visible points} & \textbf{$\alpha_i$} & \textbf{$\beta_i$} \\
        \midrule
        $\bQ\left(\sqrt2\right)$ & $\checkmark$ & $(1)$ & 1 & $\sqrt2 - 1$ \\
        \midrule
        $\bQ\left(\sqrt3\right)$ & $\times$ & $(1), (\sqrt{3} + 1)$ & $1, \frac{1}{2}$ & $\frac{1}{2}\sqrt3 - \frac{1}{2},\frac{1}{2}\sqrt3 - \frac{1}{2}$ \\
        \midrule
        $\bQ\left(\sqrt5\right)$ & $\checkmark$ & $(1)$ & 1 & $\frac{1}{2}\sqrt5 - \frac{1}{2}$ \\
        \midrule
        $\bQ\left(\sqrt6\right)$& $\times$ & $(1), (\sqrt{6} - 2)$ & $1, \frac{1}{2}$ & $\frac{1}{2}\sqrt6 - 1,\frac{1}{2}\sqrt6 - 1$ \\
        \midrule
        $\bQ\left(\sqrt7\right)$ & $\times$ &
        $\begin{aligned}
        &(1), (\sqrt{7} + 3), \\
        &(\sqrt{7} + 2), (\sqrt{7} - 2)
        \end{aligned}$
        & $1, \frac{1}{2}, \frac{1}{3}, \frac{1}{3}$
        & $\begin{aligned}
        &\frac{1}{3} \sqrt7 - \frac{2}{3}, \frac{1}{6}\sqrt7 - \frac{1}{6}, \\
        &\frac{1}{3} \sqrt7 - \frac{2}{3}, \frac{1}{6}\sqrt7 - \frac{1}{6}
        \end{aligned}$ \\
        \midrule
        $\bQ\left(\sqrt{10}\right)$ & $\times$ &
        $\begin{aligned}
        &(1), (2,\sqrt{10}), \\
        &(3, \sqrt{10} + 1), \\
        &(3, \sqrt{10} + 2)
        \end{aligned}$
        & $1, \frac{1}{2}, \frac{1}{3}, \frac{1}{3}$ &
        $\begin{aligned}
        &\sqrt{10} - 3, \frac{1}{6}\sqrt{10} - \frac{1}{3}, \\
        &\frac{1}{9}\sqrt{10} - \frac{1}{9}, \\
        &\frac{1}{6}\sqrt{10} - \frac{1}{3}
        \end{aligned}$\\
        \bottomrule
    \end{tabular}
\end{table}
\end{rmk}
\section{Two windows, with the same volume, but different density for the visible points} \label{sec:basic_rings_different_density_example}
The lattice $\cL$ from \eqref{eq:multi_dimension_number_field_lattice} gives cut and project sets with various interesting properties, such as the Amman-Beenker set in the case $K = \bQ(\sqrt{2})$ with a certain window. The question of visibility in those cut and projects is especially interesting, as the blockings in $\Lambda(\cW,\cL)$ are caused by multiplying by elements of $K$, thus understanding the visible points comes into contact with factorization in $\cO_K$.

In this section, take $K =\bQ(\sqrt{2})$ (any other real quadratic number field which satisfies H-condition works as well), and any $d \ge 2$, and we take the lattice $\cL$ from \eqref{eq:multi_dimension_number_field_lattice}. We will find two windows $\cW_1,\cW_2 \subseteq \bR^d$, such that $\vol(\cW_1) = \vol(\cW_2)$, while $\theta(\Lambda_{vis}(\cW_1,\cL)) \ne \theta(\Lambda_{vis}(\cW_2,\cL))$.
This shows the visible points might have strange properties which affect the density, and depend in weird ways on the window, and might not only depend on its volume. So just the result of proving the density always exists, from \cite[Theorem 1]{marklof2014visibilitydirectionsquasicrystals}, is surprising on its own.

Denote the closed ball of radius $R > 0$ centered at the origin in $\bR^d$ by $B_R^d$, and the closed unit ball $B^d = B_1^d$. Set $c = \sqrt{2}$ (For other $K$, $c > 1$ is a parameter which depends only on $K$). Now take $R_1,R_2,r > 0$, with $R_2 > R_1$, and $R_2^d - R_1^d = r^d$, and $\frac{R_2}{R_1} < c$, and choose the two windows:
\[
\cW_1 = B_{R_2}^d \setminus B_{R_1}^d, \ \cW_2 = B_r^d.
\]
\begin{thm} \label{thm:same_volume_different_density}
    We have $\vol(\cW_1) = \vol(\cW_2)$, but:
    \[
    \theta(\Lambda_{vis}(\cW_1,\cL)) > \theta(\Lambda_{vis}(\cW_2,\cL)).
    \]
\end{thm}
\begin{rmk}
    Note $\theta(\Lambda_{vis}(\cW_1,\cL)), \theta(\Lambda_{vis}(\cW_2,\cL)),$ both exist and don't depend on the averaging set $D$, by \cite[Theorem 1]{marklof2014visibilitydirectionsquasicrystals}.
\end{rmk}
The volume equality is easy:
\[
\vol(\cW_1) = R_2^d \vol(B^d) - R_1^d \vol(B^d) = r^d \vol(B^d) = \vol(\cW_2).
\]
By \cite[Theorem 4.9]{hammarhjelm2019densityminimalgapvisible} (the proof there generalizes to any $d \ge 2$), because $\cW_2$ is centrally symmetric, convex, and bounded, we have: 
\[
\theta \left( \Lambda_{vis}(\cW_2,\cL)\right) = \left( 1 - \frac{1}{\lambda^d} \right) \frac{\vol(\cW_2)}{\zeta_K(d)\covol(\cL)}.
\]

Where $\lambda = 1 + \sqrt2 > 1$ is the fundamental unit as before. As $\cW_1$ is not convex (more importantly, not star-shaped with respect to the origin), the previous theorem does not apply to it. 
\begin{lem}
    We have:
    \[\Lambda_{prim}(\cW_1,\cL,(1)) \subseteq \Lambda_{vis}(\cW_1,\cL).\]
\end{lem}

\begin{proof}
    Let $\bm{x} \in \Lambda_{prim}(\cW_1, \cL, (1))$, and assume by way of contradiction $\bm{y} = t \bm{x} \in \Lambda(\cW_1,\cL)$ blocks $\bm{x}$. We must have $0 < t < 1$, and also $t \in K$. Assume $t = \frac{\alpha}{\beta}$, where $\alpha,\beta \in \cO_K$, and $\gcd(\alpha,\beta) = (1)$ (we can assume this as $\cO_K$ is a unique factorization domain). For all $1 \le i \le d$, we have $\beta y_i = \alpha x_i$, and as $\gcd(\alpha,\beta) = (1)$, we must have $\beta | x_i$. So because $\gcd(x_1,...,x_d) = (1)$ (by the definition of the primitive points), we must have that $\beta$ is a unit. So from the start, $t \in \cO_K$.
    Now notice because both $\bm{x},\bm{y} \in \Lambda(\cW_1,\cL)$, we must have $\sigma(\bm{x}),\sigma(\bm{y}) \in \cW_1$, and $\sigma(\bm{y}) = \sigma(t) \sigma(\bm{x})$, thus by the definition of $\cW_1$:
    \[
    \left\{
    \begin{aligned}
        &R_1 < \|\sigma(\bm{x})\| \le R_2 \\
        &R_1 < |\sigma(t)| \cdot \|\sigma(\bm{x})\| \le R_2
    \end{aligned}
    \right.
    \Rightarrow \frac{R_1}{R_2} < |\sigma(t)| < \frac{R_2}{R_1} \Rightarrow \frac{1}{c} < |\sigma(t)| < c.
    \]
    We also have $0 < t < 1$, and $t \in \cO_K$. It is now a simple case check such $t$ doesn't exist. As $\cO_{\bQ\left(\sqrt{2}\right)} = \bZ[\sqrt{2}]$, we can write $t = a + b\sqrt{2}$ for $a,b \in \bZ$, and we have:
    \[
    \left\{
    \begin{aligned}
        &0 < a + b\sqrt{2} < 1, \\
        &-\sqrt{2} < a-b\sqrt{2} < \sqrt{2}.
    \end{aligned}
    \right.
    \Rightarrow
    \left\{
    \begin{aligned}
        &-\frac{\sqrt{2}}{2} < a < \frac{\sqrt{2}}{2} + \frac{1}{2}, \\
        &-\frac{1}{2} < b < \frac{1}{2} + \frac{1}{2 \sqrt{2}}.
    \end{aligned}
    \right.
    \]
    By the second condition, $b = 0$, but if $b = 0$ we can't have $0 < a + b\sqrt{2} < 1$. So we don't have any such $t$. 
\end{proof}
\begin{rmk}
    For other $K$ with the H-condition, we can find small enough $c > 1$ such that there is no $t \in \cO_K$ with $0 < t < 1$ and $\frac{1}{c} < |\sigma(t)| < c$.
\end{rmk}
\begin{rmk}
    Notice the difference the lemma gives us- all of the primitive points are also visible, unlike with star-shaped windows, where we must remove a portion for which $\sigma(\bm{x}) \in \frac{1}{\lambda}\cW$ as in \Cref{prop:h_cond_visible_description}.
\end{rmk}
Now we have $ \Lambda_{prim}(B_{R_1}^d,\cL,(1)) \subseteq \Lambda_{prim}(B_{R_2}^d,\cL,(1))$, and also:
\[
\Lambda_{vis}(\cW_1,\cL) \supseteq \Lambda_{prim}(\cW_1,\cL,(1)) = \Lambda_{prim}\left(B_{R_2}^d,\cL,(1)\right) \setminus \Lambda_{prim}\left(B_{R_1}^d,\cL,(1)\right) .
\]
\begin{proof} [Proof of \Cref{thm:same_volume_different_density}]
    Using the calculation of the density of the primitive points (\Cref{lem:primitive_point_counting}, also proved in \cite{hammarhjelm2019densityminimalgapvisible} for those cases), now for the balls which are convex sets, we get:
    \[
        \begin{aligned}
            &\theta\left(\Lambda_{vis}(\cW_1,\cL)\right) \ge \theta\left(\Lambda_{prim}(\cW_1,\cL,(1))\right) = \\
            &= \theta\left(\Lambda_{prim}\left(B_{R_2}^d,\cL,(1)\right)\right) - \theta\left(\Lambda_{prim}\left(B_{R_1}^d,\cL,(1)\right) \right) = \\
            &= \frac{\vol\left(B_{R_2}^d\right)}{\zeta_K(d) \covol(\cL)} - \frac{\vol\left(B_{R_1}^d\right)}{\zeta_K(d) \covol(\cL)} = \frac{\vol(\cW_1)}{\zeta_K(d)\covol(\cL)} > \\ 
            &> \left( 1 - \frac{1}{\lambda^d} \right) \frac{\vol(\cW_2)}{\zeta_K(d)\covol(\cL)} = \theta\left(\Lambda_{vis}(\cW_2,\cL)\right).
        \end{aligned}
    \]
    Thus the density of the visible points is different.
\end{proof}
\section{Range of densities of the visible points when changing the window, with a random lattice} \label{sec:density_range_random_lattice}

\subsection{Notations}
In this section we will go back to the more general setting, where $m,d$ can be different, and our windows and averaging sets will not necessarily be nice.

Denote by $\mu$ the Haar-measure on $GL_n(\bR)$, and we will look at the lattice defined by $\cL_g = g \bZ^n$, where $g$ is chosen to be outside a set of $\mu$ measure 0.

First, let us define the set of $\mu$-full measure for which we prove \Cref{thm:random_density_range}.
\begin{defn} \label{def:full_measure_generic_set}
    For every $\bm{u},\bm{v} \in \bZ^n$, define:
    \[
    \sB(\bm{u},\bm{v}) = \{g \in GL_n(\bR) | \pi_{\rm phys}(g \bm{u}) , \pi_{\rm phys}(g \bm{v}) \textit{ are linearly dependent} \}.
    \]
    Also define:
    \[
    \sB = \bigcup_{\substack{\bm{u},\bm{v} \in \bZ^n \\ \textit{linearly independent}}} \sB(\bm{u},\bm{v}).
    \]
    We take:
    \[
    \sA = \left\{ g \in GL_n(\bR) : g\notin\sB,\overline{\pi_{\rm int}(\cL_g)} =\bR^m, \pi_{\rm phys}|_{\cL_g}\textit{ is injective}\right\},
    \]
    i.e. $\sA$ is all $g \in GL_n(\bR) \setminus \sB$ for which $\cL_g$ gives irreducible cut and project sets.
\end{defn}
\begin{prop}
    We have:
    \[\mu\left(GL_n(\bR) \setminus \sA\right) = 0.\]
\end{prop}
\begin{proof}
    For each $\bm{u},\bm{v}$, the set $\sB(\bm{u},\bm{v})$ is an algebraic subvariety of $GL_n(\bR)$, and as $d \ge 2$, if $\bm{u},\bm{v}$ are independent it is a proper subvariety. It is known that proper algebraic subvarieties of $GL_n(\bR)$ are of $\mu$-measure 0. Thus $\mu(\sB(\bm{u},\bm{v})) = 0$ for independent $\bm{u},\bm{v} \in \bZ^n$. Thus as a countable union of measure 0 sets, also $\mu(\sB) = 0$.
    
    As $\sB$ is of zero measure, and for $\mu \textit{-a.e. } g$ we have that $\cL_g$ gives irreducible cut and project sets, $\sA$ is a set of full measure.
\end{proof}
We will prove \Cref{thm:random_density_range,thm:random_density_formula} for all $g \in \sA$.

\subsection{Visible points and density formula}
For any $g \in \sA$, and any window $\cW$, we will find a description of the visible points, and from that, we will prove \Cref{thm:random_density_formula}.

We say a point $\bm{y} \ne 0$ blocks a point $\bm{x} \ne 0$, if there exists $t \in (0,1)$ such that $\bm{y} = t \bm{x}$. If for a pair $\bm{x},\bm{y}$ such $t$ exists, it is unique.

\begin{lem} \label{lem:pairwise_blocking_condition_random}
    For all $g \in \sA$, we have that for any $\bm{x},\bm{y} \in \bZ^n \setminus \{ 0 \}$:
    \[
    \bm{y} \textit{ blocks } \bm{x} \iff \pi_{\rm phys}(g \bm{y}) \textit{ blocks } \pi_{\rm phys}(g \bm{x}).
    \]
\end{lem}
\begin{rmk}
    The argument for this lemma was already in \cite[proposition 2.2]{gringlaz2025errorboundsvisiblepoints}. We provide a proof for completeness.
\end{rmk}
\begin{proof}
    Let $g \in \sA$. First, note that for any $\bm{x} \in \bZ^n \setminus {0}$, we have $\pi_{\rm phys}(g \bm{x}) \ne 0$.
    
    Now the $\Rightarrow$ direction is obvious, as if $\bm{y} = t \bm{x}$, we also have $\pi_{\rm phys}(g \bm{y}) = t \pi_{\rm phys}(g \bm{x})$.
    
    For the $\Leftarrow$ direction, assume $\bm{y}$ doesn't block $\bm{x}$.
    If they are linearly independent, we have $g \notin \sB(\bm{x},\bm{y})$ by \Cref{def:full_measure_generic_set}, thus $\pi_{\rm phys}(g \bm{y}), \pi_{\rm phys}(g \bm{x})$ are linearly independent, and can't block each other.
    If $\bm{x},\bm{y}$ are linearly dependent, then $\bm{y} = t \bm{x}$ for $t \notin (0,1)$, thus $\pi_{\rm phys}(g \bm{y}) = t \pi_{\rm phys}(g \bm{x})$ for $t \notin (0,1)$, so the blocking doesn't occur.
\end{proof}

As we understand visibility in $\bZ^n$, using the lemma, we can use it to understand visibility in $\Lambda(\cW, \cL_g)$, as can be seen in the following proposition:
\begin{prop} \label{prop:random_visible_classification}
    For any $k \in \bN$, denote $\bZ_{k}^n = \{ \bm{x} \in \bZ^n \ | \ \gcd(x_1,...,x_n)=k\}$. We have:
    \[
        \Lambda_{vis}(\cW, \cL_g) = \biguplus_{k=1}^{\infty}{\Lambda\left( \cW \setminus \left( \bigcup_{i=1}^{k-1}\frac{k}{i}\cW \right), g \bZ_{k}^n \right)}.
    \]
    And the union is disjoint.
\end{prop}
\begin{proof}
     We have $\bZ_\star^n = \biguplus_{k=1}^{\infty}\bZ_k^n$. For $\bm{x} \in \bZ_k^n$, the points which block it in $\bZ^n$ are $\{\frac{i}{k} \bm{x}\}_{i=1}^{k-1}$, thus by \Cref{lem:pairwise_blocking_condition_random}, the points which block $\pi_{\rm phys}(g\bm{x})$ in $\pi_{\rm phys}(\cL_g)$ are $\{\frac{i}{k} \pi_{\rm phys}(g\bm{x})\}_{i=1}^{k-1}$. So $\pi_{\rm phys}(g\bm{x}) \in \Lambda_{vis}(\cW, \cL_g)$ if and only if $\frac{i}{k} \pi_{\rm int}(\bm{x}) \notin \cW$ for $i=1,...,k-1$, or equivalently, $\pi_{\rm int}(\bm{x}) \notin \frac{k}{i} \cW$ for those $i$. So we have for all windows $\cW$:
    \[
        \begin{aligned}
            &\Lambda_{vis}(\cW, \cL_g) = \biguplus_{k=1}^{\infty}{\left(\Lambda_{vis}(\cW, g \bZ^n) \cap \Lambda\left( \cW, g \bZ_{k}^n \right)\right)} = \\
            &= \biguplus_{k=1}^{\infty}{\Lambda\left( \cW \setminus \left( \bigcup_{i=1}^{k-1}\frac{k}{i}\cW \right), g \bZ_{k}^n \right)}.
        \end{aligned}
    \]
\end{proof}

\begin{lem} \label{lem:bound_point_count_in_multiples}
    For any $g \in \sA$, and bounded sets $\cW \subset \bR^m, D \subset \bR^d$, there exists $C > 0$, such that for any $i \in \bN$, and $T > 1$ we have:
    \[
        \#\left( \Lambda_\star(\cW, ig\bZ^n) \cap TD \right) \le \frac{CT^d}{i^d}.
    \]
\end{lem}
\begin{proof}
     Inscribe $\cW, D$ inside constant boxes $\cW', D'$ such that $\bm{0} \in \cW' \times D'$, and take a small enough constant box $B \subset \bR^n$ whose diameter is less than the shortest vector in $g \bZ^n$. We can take $C > 0$ large enough so for any $i \in \bN$, $T > 1$, we can cover $\cW \times TD \subseteq \cW' \times TD'$ with at most $\left\lceil C \left(\frac{T}{i}\right)^d\right\rceil$ copies of $iB$.
     Each copy contains at most one vector from $ig\bZ^n$, and note one of the copies contains the zero vector which we don't count. Thus:
     \[
     \#\left( \Lambda_\star(\cW, ig\bZ^n) \cap TD \right) \le \#((\cW' \times TD') \cap ig\bZ^n_\star) \le \left\lceil \frac{CT^d}{i^d} \right\rceil  - 1 < \frac{CT^d}{i^d}.
     \]
\end{proof}

\begin{lem}
    For any bounded set $\cW \subseteq \bR^m$ with boundary of measure zero, any $g \in \sA$, and any $k \in \bN$ the set $\Lambda(\cW,g\bZ_k^n)$ has density $\frac{\vol(\cW)}{k^n \zeta(n) \covol(g \bZ^n)}$.
\end{lem}

\begin{proof}
    Firstly, we can assume $k = 1$, as $g\bZ_k^n = (kg) \bZ_1^n$, and we can take $kg$ instead of $g$. Also assume $\cW$ has non-empty interior, we will explain the other case later. This is close to well-known propositions (see \cite[Proposition 6]{baake2000diffraction}), where the number of points in $g\bZ_1^n$ is counted along families of the form $TB$ for $B \subseteq \bR^n$ Jordan-measurable, and growing $T\rightarrow\infty$. But here we count the number of points in $g\bZ_1^n$, along the family $\cW \times TD$, which is not a typical ordered family to average along, thus we need to use the irreducibility properties of $\cL_g$ to get the same expected density.

    Let $D \subseteq \bR^d$ be an averaging set, and denote $\mu: \bN \rightarrow \bZ$ the standard M\"obius function:
    \[
    \mu(i) = \begin{cases} 
    1 & \text{if } i = 1, \\
    (-1)^r & \text{if } i = p_1 \cdots p_r \text{ for distinct primes } p_1,...,p_r, \\
    0 & \text{if } p^2|i \text{ for some prime } p.
    \end{cases}
    \]
    As $\pi_{\rm phys} \circ g : \bZ^n \rightarrow \bR^d$ is injective, and $\Lambda(\cW,g\bZ_1^n) \cap TD$ is finite, we have by inclusion exclusion:
    \[
    \# \left( \Lambda(\cW,g\bZ_1^n) \cap TD \right) = \sum_{i=1}^{\infty}{\left( \mu(i)\cdot\#\left( \Lambda_\star(\cW,ig\bZ^n) \cap TD \right) \right)}.
    \]
    As $\cL_g$ gives irreducible cut and project sets, for any $i \in \bN$ the map $\pi_{\rm phys}|_{i \cL_g}$ is injective as $i \cL_g \subseteq \cL_g$, and $\pi_{\rm int}\left(i \cL_g\right) = i \pi_{\rm int}\left(\cL_g\right)$ is dense in $\bR^m$. So the set $\Lambda(\cW,ig\bZ^n)$ is an irreducible cut and project set, so it has density $\frac{\vol(\cW)}{\covol(ig \bZ^n)} = \frac{\vol(\cW)}{\covol(g \bZ^n)} \cdot \frac{1}{i^n}$.

    We have that $\frac{1}{\zeta(n)} = \sum_{i=1}^{\infty}{\frac{\mu(i)}{i^n}}$, for any $N_0 \in \bN$, as $\cW,D$ are bounded, we can use \Cref{lem:bound_point_count_in_multiples} in order to bound the tail in the sum:
    \[
    \begin{aligned}
        & \left| \frac{\#\left( \Lambda(\cW,g\bZ_1^n) \cap TD \right)}{\vol(TD)} -  \frac{\vol(\cW)}{\zeta(n) \covol(g \bZ^n)} \right| \le \\
        & \le \sum_{i=1}^{\infty}{\left| \frac{\#\left( \Lambda_\star(\cW,ig\bZ^n) \cap TD \right)}{\vol(TD)} - \frac{\vol(\cW)}{i^n \covol(g \bZ^n)}\right|} \le \\
        &\sum_{i=1}^{N_0}{\left| \frac{\#\left( \Lambda_\star(\cW,ig\bZ^n) \cap TD \right)}{\vol(TD)} - \frac{\vol(\cW)}{i^n \covol(g \bZ^n)}\right|} + \sum_{i=N_0+1}^{\infty}{\left(\frac{CT^d}{i^d \vol(TD)} + \frac{\vol(\cW)}{i^n \covol(g \bZ^n)}\right)} \\
        &\le \sum_{i=1}^{N_0}{\left| \frac{\#\left( \Lambda_\star(\cW,ig\bZ^n) \cap TD \right)}{\vol(TD)} - \frac{\vol(\cW)}{i^n \covol(g \bZ^n)}\right|} + C' \sum_{i=N_0+1}^{\infty}{\frac{1}{i^d}}.
    \end{aligned}
    \]
    Where $C'$ doesn't depend on $T, N_0$, so as we have $\theta\left(\Lambda_\star(\cW,ig\bZ^n) \right) = \frac{\vol(\cW)}{i^n \covol(g \bZ^n)}$, when taking $T \rightarrow \infty$, the first $N_0$ summands go to zero, so:
    \[
    \limsup_{T\rightarrow\infty} \left| \frac{ \# \left( \Lambda(\cW,g\bZ_1^n) \cap TD \right)}{\vol(TD)} - \frac{\vol(\cW)}{\zeta(n) \covol(g \bZ^n)}\right| \le C' \sum_{i=N_0+1}^{\infty}{\frac{1}{i^d}}\overset{N_0 \rightarrow \infty}{\rightarrow} 0.
    \]
    So indeed $\theta(\Lambda(\cW,g\bZ_1^n)) = \frac{\vol(\cW)}{\zeta(n) \covol(g \bZ^n)}$.

    If $\cW$ has empty interior, we have $\bar{\cW} = \partial \cW$ is a compact set of measure 0, thus we can approximate it from above for any $\eps > 0$ by a finite union of rectangles $\cW_\eps$ who contains $\bar{\cW}$, and $\vol(\cW_\eps) < \eps$, so for any averaging set $D$:
    \[
    \limsup_{T\rightarrow\infty} \frac{ \# \left( \Lambda(\cW,g\bZ_1^n) \cap TD \right)}{\vol(TD)} \le \theta(\Lambda(\cW_\eps,g\bZ_1^n)) < \frac{\eps}{\zeta(n) \covol(g \bZ^n)} \overset{\eps \rightarrow 0}{\rightarrow} 0
    \]
    So we also get $\theta(\Lambda(\cW,g\bZ_1^n)) = 0$.
\end{proof}

\begin{proof} [Proof of \Cref{thm:random_density_formula}]
    For $k \in \bN$, denote $\cP_k = \Lambda\left( \cW \setminus \left( \bigcup_{i=1}^{k-1}\frac{k}{i}\cW \right), g \bZ_{k}^n \right)$, and note they are all disjoint, and $\Lambda_{vis}(\cW, \cL_g) = \biguplus_{k=1}^{\infty}{\cP_k}$. As $\cW \setminus \left( \bigcup_{i=1}^{k-1}\frac{k}{i}\cW \right)$ is bounded with boundary measure 0, by the previous lemma, we have $\theta(\cP_k) = \frac{\vol\left( \cW \setminus \left( \bigcup_{i=1}^{k-1}\frac{k}{i}\cW \right) \right)}{k^n \zeta(n) \covol(g \bZ^n)}$. Let $D \subset \bR^d$ be some averaging set. We can see that $\cP_k \subseteq  \Lambda\left(\cW, k g\bZ^n \right)$, thus by \Cref{lem:bound_point_count_in_multiples}, we have that $ \#(\cP_k \cap TD) \le \frac{CT^d}{k^d}$, with $C>0$ independent of $T,k$. Thus for any $N_0 \in \bN$, we have:
    \[
    \begin{aligned}
        &\left| \frac{\#(\Lambda_{vis}(\cW, \cL_g) \cap TD)}{\vol(TD)} - \frac{1}{\zeta(n) \covol(\cL_g)}\sum_{k=1}^{\infty}{\frac{1}{k^n} \vol\left( \cW \setminus \left( \bigcup_{i=1}^{k-1}\frac{k}{i}\cW \right)\right)} \right| \le \\
        & \le \sum_{k=1}^{\infty}{\left| \frac{\#\left( \cP_k \cap TD \right)}{\vol(TD)} - \frac{\vol\left( \cW \setminus \left( \bigcup_{i=1}^{k-1}\frac{k}{i}\cW \right)\right)}{\zeta(n) k^n \covol(\cL_g)}\right|} \le \\
        &\sum_{k=1}^{N_0}{\left| \frac{\#\left( \cP_k \cap TD \right)}{\vol(TD)} - \theta(\cP_k)\right|} + \sum_{k=N_0+1}^{\infty}{\left(\frac{CT^d}{k^d \vol(TD)} + \frac{\vol(\cW)}{k^n \covol(\cL_g)}\right)} \\
        &\le \sum_{k=1}^{N_0}{\left| \frac{\#\left( \cP_k \cap TD \right)}{\vol(TD)} - \theta(\cP_k)\right|}  + C' \sum_{k=N_0+1}^{\infty}{\frac{1}{k^d}}.
    \end{aligned}
    \]
    Where $C'$ doesn't depend on $T, N_0$. When taking $T \rightarrow \infty$, the first $N_0$ summands go to zero, so:
    \[
    \begin{aligned}
        &\limsup_{T\rightarrow\infty} \left| \frac{\#(\Lambda_{vis}(\cW, \cL_g) \cap TD)}{\vol(TD)} - \frac{1}{\zeta(n) \covol(\cL_g)}\sum_{k=1}^{\infty}{\frac{1}{k^n} \vol\left( \cW \setminus \left( \bigcup_{i=1}^{k-1}\frac{k}{i}\cW \right)\right)} \right| \le\\
        & \le C' \sum_{k=N_0+1}^{\infty}{\frac{1}{k^d}}\overset{N_0 \rightarrow \infty}{\rightarrow} 0.
    \end{aligned}
    \]
    So we get $\Lambda_{vis}(\cW, \cL_g)$ has the density we claimed.
\end{proof}

\subsection{Possible range of the density, \Cref{thm:random_density_range}}
Now we are ready to prove the bounds in \Cref{thm:random_density_range} (1):
\begin{proof} [Proof of \Cref{thm:random_density_range} (1)]
    Consider $g \in \sA$. By considering only $k = 1$ in \Cref{thm:random_density_formula}, we get:
    \[
    \theta(\Lambda_{vis}(\cW, \cL_g)) \ge \frac{1}{\zeta(n)} \theta(\Lambda(\cW, \cL)).
    \]
    For the upper bound, note for any $k \in \bN$ we can bound $\vol\left(\cW\setminus \left( \bigcup_{i=1}^{k-1}\frac{k}{i}\cW \right)\right) \le \vol(\cW)$. But also, because $\cW$ has non-empty interior, there exists a rational number $\frac{p}{q} > 1$ such that $\vol\left( \cW \cap \frac{p}{q} \cW \right) > 0$, so for $k = p$, we have
    \[
    \vol\left(\cW\setminus \left( \bigcup_{i=1}^{k-1}\frac{k}{i}\cW \right)\right) \le \vol\left(\cW \setminus \frac{p}{q} \cW\right) < \vol(\cW).
    \]
    So we have strict inequality:
    \[
    \begin{aligned}
        &\theta(\Lambda_{vis}(\cW, \cL_g)) = \frac{1}{\zeta(n) \covol(\cL_g)}\sum_{k=1}^{\infty}{\frac{1}{k^n} \vol\left( \cW \setminus \left( \bigcup_{i=1}^{k-1}\frac{k}{i}\cW \right)\right)} < \\
        &< \frac{1}{\zeta(n) \covol(\cL_g)}\sum_{k=1}^{\infty}{\frac{1}{k^n} \vol(\cW)} = \theta(\Lambda(\cW, \cL)).
    \end{aligned}
    \]
\end{proof}

Now for the examples, denote $B_R^m \subset \bR^m$ the closed ball centered at the origin of radius $R \ge 0$, and $B^m=B_1^m$. Fix $V>0$ for the rest of the proof, for a parameter $c \in (1, \infty]$, we will consider:
\[
\cW_c = B_{R_2}^m \setminus B_{R_1}^m,
\]
where $R_1 \ge 0, R_2 > 0$ are defined uniquely by $(R_2^m-R_1^m)\vol(B^m) = V$, and $\frac{R_2}{R_1} = c$, and in the case $c = \infty$, we take $R_1 = 0$ and $R_2^m \vol(B^m) = V$. We have $\vol(\cW_c) = V$, and we will show $\left\{ \theta(\Lambda_{vis}(\cW_c, \cL_g)) : c \in (1,\infty] \right\}$ gives us the entire range from \Cref{thm:random_density_range}. 

\begin{lem} \label{lem:w_c_remove_invisible_volume}
    Let $\cW_c$ be as above. If $c < \infty$, then for $k \in \bN$, we have:
    \[
    \vol\left( \cW_c \setminus \left( \bigcup_{i=1}^{k-1}\frac{k}{i} \cW_c \right)\right) =
    \begin{cases} 
        \frac{\left( \frac{k}{k-1} \right)^m - 1}{c^m - 1} V & \text{if } k \ge 1+\frac{1}{c-1}, \\
        V & \text{if } k \le 1+\frac{1}{c-1}.
    \end{cases}
    \]
    And if $c = \infty$:
    \[
    \vol\left( \cW_c \setminus \left( \bigcup_{i=1}^{k-1}\frac{k}{i} \cW_c \right)\right) =
    \begin{cases} 
        0 & \text{if } k > 1, \\
        V & \text{if } k = 1.
    \end{cases}
    \]
\end{lem}
\begin{proof}
    If $c = \infty$, we have $\cW_c \subseteq t \cW_c$ for any $t \ge 1$, and the result follows.
    
    If $c < \infty$, from the shape of $\cW_c$, we have for any $i=1,..,k-1$ that:
    \[
    \frac{k}{i} \cW_c \cap B_{R_2}^m \subseteq \frac{k}{k-1} \cW_c \cap B_{R_2}^m.
    \]
    Thus also:
    \[
    \frac{k}{i} \cW_c \cap \cW_c \subseteq \frac{k}{k-1} \cW_c \cap \cW_c.
    \]
    So:
    \[
    \cW_c \setminus \left( \bigcup_{i=1}^{k-1}\frac{k}{i} \cW_c \right) = \cW_c \setminus \frac{k}{k-1} \cW_c.
    \]
    
    We have $k \le 1 + \frac{1}{c-1} \iff \frac{k}{k-1} \ge c$, so in this case, 
    \[\frac{k}{k-1} \cW_c \subset \frac{k}{k-1} (\bR^m \setminus B_{R_1}^m) \subseteq c (\bR^m \setminus B_{R_1}^m) \subseteq \bR^m \setminus \cW_c.
    \]
    So:
    \[
    \vol\left( \cW_c \setminus \left( \bigcup_{i=1}^{k-1}\frac{k}{i} \cW_c \right)\right) = \vol\left( \cW_c \setminus \frac{k}{k-1} \cW_c \right) = \vol(\cW_c) = V.
    \]

    In the other case, we have $k \ge 1 + \frac{1}{c-1} \iff \frac{k}{k-1} \le c$, and noting $(c^m - 1) R_1^m \vol(B^m) = V$, we get:
    \[
    \begin{aligned}
        &\vol\left( \cW_c \setminus \left( \bigcup_{i=1}^{k-1}\frac{k}{i} \cW_c \right)\right) = \vol\left( \cW_c \setminus \frac{k}{k-1} \cW_c \right) = \\
        &= \vol\left(B_{\frac{k}{k-1} R_1}^m \setminus B_{R_1}^m \right) = \left( \left( \frac{k}{k-1} \right)^m - 1\right) R_1^m \vol(B^m) = \frac{\left( \frac{k}{k-1} \right)^m - 1}{c^m - 1} V.
    \end{aligned}
    \]
\end{proof}

\begin{proof} [Proof of \Cref{thm:random_density_range} (2)]
    For $k \in \bN$, define $f_k: (1,\infty] \rightarrow \bR_{\ge0}$ by:
    \[
    f_k(c) = \vol\left( \cW_c \setminus \left( \bigcup_{i=1}^{k-1}\frac{k}{i} \cW_c \right)\right).
    \]
    By \Cref{lem:w_c_remove_invisible_volume}, $f_k$ is decreasing and continuous, including at infinity, $\lim_{c \rightarrow \infty}f_k(c) = f_k(\infty)$. Looking at the sum:
    \[
    F(c) := \sum_{k=1}^\infty{\frac{1}{k^n}f_k(c)}  = \theta\left(\Lambda_{vis}\left(\cW_c, \cL_g\right)\right) \cdot\zeta(n) \covol(\cL_g).
    \]
    We claim $F|_{(1,\infty)}$ is continuous. For $N_0 \in \bN$, when restricting $F|_{\left( \frac{N_0}{N_0-1},\infty \right)}$, as $c > \frac{N_0}{N_0-1}$, we have $N_0 > 1 + \frac{1}{c-1}$, and the sum can also be written as:
    \[
    F(c) = \sum_{k = 1}^{N_0-1}{\frac{1}{k^n}f_k(c)} + \sum_{k=N_0}^{\infty}{\frac{1}{k^n}\left( \frac{\left( \frac{k}{k-1} \right)^m - 1}{c^m - 1} V\right)} = \sum_{k = 1}^{N_0-1}{\frac{1}{k^n}f_k(c)} + \frac{1}{c^m - 1} \cdot C_{N_0}V.
    \]
    Where $C_{N_0}$ is a constant independent of $c$. This is a sum of finitely many continuous functions, thus it is continuous. So $F|_{\left( \frac{N_0}{N_0-1},\infty \right)}$ is continuous for any $N_0 \in \bN$, so $F|_{(1,\infty)}$ is continuous.
    
    Also note as all of the $f_k$ are non-negative and decreasing, we also have $\lim_{c \rightarrow \infty}{F(c)} = F(\infty)$ and $\lim_{c \rightarrow 1^+}{F(c)} = \sum_{k=1}^\infty{\frac{1}{k^n}\lim_{c \rightarrow 1^+}{f_k(c)}}$.

    Thus as $F$ is continuous, it achieves every value between $\lim_{c \rightarrow 1^+}{F(c)}$ and $F(\infty)$, which are:
    \[
    \begin{aligned}
        &\lim_{c \rightarrow 1^+}{F(c)} = \sum_{k=1}^\infty{\frac{1}{k^n}\lim_{c \rightarrow 1^+}{f_k(c)}} =\sum_{k=1}^{\infty}{\frac{1}{k^n}V} = \zeta(n)V, \\
        &F(\infty) = \sum_{k=1}^{\infty}{\frac{1}{k^n}f_k(\infty)} = \frac{1}{1^n} \cdot V + \sum_{k=2}^{\infty}{\frac{1}{k^n} \cdot 0} = V.
    \end{aligned}
    \]
    So for any $q \in \left[ \frac{1}{\zeta(n)}\frac{V}{\covol(\cL_g)}, \frac{V}{\covol(\cL_g)} \right)$, there is a $c \in (1, \infty]$ such that $F(c) = q \cdot \zeta(n) \covol\left(\cL_g\right)$, and  $\theta\left( \Lambda_{vis}(\cW_c, \cL_g) \right) = q$.
\end{proof}
\section{Bounds on the distance of holes in the set of visible points for the H-condition cases} \label{sec:holes}
In this section we will mostly go back to the same notations as \Cref{sec:non_h_density}. Let $K$ be a quadratic real number field satisfying H-condition, take some $d \ge 2$, and define $\cL$ by \eqref{eq:multi_dimension_number_field_lattice}. Let $\cW \subseteq \bR^d$ be a nice window for $d \ge 2$, and denote $D_0 = [0,1]^d$, and we will give bounds on finding $\bm{x} \in \bR^d$ with $(\bm{x} +TD_0) \bigcap \Lambda_{vis}(\cW,\cL) = \es$.

\subsection{Upper bound on the distance,  \Cref{thm:large_holes_distance} (1)}
\begin{lem} \label{lem:lattice_point_in_rectangle}
    There is a constant $C > 0$ depending only on $K$, such that for all rectangles $R \subset \bR^2$ (with sides parallel to the axis), if $\vol(R) \ge C$, then $R \cap \cL_K \ne \es$, where $\cL_K = \{ (x, \sigma(x)) : x \in \cO_K \}$ is the Minkowski embedding of $\cO_K$.
\end{lem}
\begin{rmk}
    This is a well known and an interesting property of the lattices constructed by number fields, with relation to the Danzer problem, proved in greater generality in \cite{simmons2016danzer}. We provide a proof for completeness.
\end{rmk}
\begin{proof}
    If we restrict to $R$ being a square, this theorem is true for general lattices, and is easy to prove. Denote $C' > 0$ the constant for squares, and $C = \lambda^2 C'$. Take the linear transformation $a_\lambda = \text{diag}(\lambda, \sigma(\lambda))$ on $\bR^2$. It preserves $\cL_K$, and  for rectangles multiplies one of the side lengths by $\lambda$, while the other side is multiplied by $|\sigma(\lambda)| = \frac{1}{\lambda}$. Thus if $R \subset \bR^2$ is some rectangle with $\vol(R) \ge C$, there exists $k \in \bZ$ such that $a_\lambda^k R$ is a rectangle whose ratio of the longer side length to smaller one is less than $\lambda^2$.
    
    Denote $l_1 \le l_2$ the side lengths of $a_\lambda^k R$. Then we have $l_1^2 \ge \frac{1}{\lambda^2} l_1 l_2 = \frac{1}{\lambda^2}\vol(a_\lambda^kR) = \frac{1}{\lambda^2}\vol(R) \ge C'$. And $a_\lambda^k R$ contains a square of volume $l_1^2$, thus by the definition of $C'$ it contains a point $(x,\sigma(x)) \in \cL_K$, which means $(\lambda^{-k}x,\sigma(\lambda^{-k}x)) \in \cL_K \cap R$.
\end{proof}

Denote $U_\delta = [-\delta, \delta]^d$ for all $\delta > 0$. 
\begin{lem} \label{lem:approximate_volume} 
    Let $A \subset \bR^d$ be a bounded Jordan measurable set of positive volume. For every $\eps > 0$, there exist Jordan measurable sets $A^- \subset A \subset A^+$ such that for small enough $\delta > 0$ :
    \[
    \begin{aligned}
        &\vol(A^+) - \vol(A^-) \le \eps, \\
        &A^- + U_\delta \subseteq A,\  A + U_\delta \subseteq A^+.
    \end{aligned}
    \]
\end{lem}
\begin{proof}
    We can cover the boundary of $A$ by a union of finitely many axis aligned boxes with volume at most $\frac{\eps}{2 \cdot 3^d}$, denote the boxes $B_1,...,B_k$. Rescale each box from its center by multiplying each side by $3$, and denote the rescaled boxes $B_1',...,B_k'$. We have $\vol{\left(\bigcup_{i=1}^{k}B_i'\right)} \le \frac{\eps}{2}$. Define:
    \[
    A^- = A \setminus \left(\bigcup_{i=1}^{k}B_i'\right), A^+ = A \cup \left(\bigcup_{i=1}^{k}B_i'\right).
    \]
    This approximates $A$ from above and below with Jordan measurable sets $A^- \subseteq A \subseteq A^+$, we have $\vol(A^+) - \vol(A^-) \le \eps$. Taking $\delta > 0$ to be smaller than all the sides of $B_1,...,B_k$, the second part follows from the rescaling of the boxes, as now for any $\bm{x} \in \partial A$, $\bm{x} + U_\delta$ is contained in $\bigcup_{i=1}^{k}B_i'$, thus $A+U_\delta \subseteq A^+$ and $(A^c+U_\delta) \cap A^- = \emptyset$ and the lemma follows.
\end{proof}

\begin{proof} [Proof of \Cref{thm:large_holes_distance} (1)]
    
Denote $\eps_0 = \sqrt[5]{1 + \eps} - 1$. This is a constant as it depends only on $\eps$. Throughout the proof, $1 + \eps_0 > 1$ will appear as multiplicative constant multiple times.

We can use \Cref{lem:approximate_volume} to approximate $\cW$ from above with a set $A$, and $\frac{1}{\lambda} \cW$ from below with a set $B$, such that together:
\begin{equation} \label{eq:A_B_vol_bound}
\vol(A) - \vol(B) \le (1 + \eps_0) \left( 1 - \frac{1}{\lambda^d} \right)\vol(\cW),
\end{equation}
\begin{equation} \label{eq:A_and_B_have_small_neighborhoods}
\cW + U_\delta \subseteq A,\  B + U_\delta \subseteq \frac{1}{\lambda} \cW,
\end{equation}
for some $1 > \delta > 0$ constant (doesn't depend on $T$). And clearly $B \subseteq \frac{1}{\lambda} \cW \subseteq \cW \subseteq A$. Also denote $G = A \setminus B$, and by \eqref{eq:A_B_vol_bound} and \eqref{eq:A_and_B_have_small_neighborhoods}, we have:
\begin{equation} \label{eq:G_vol_bound}
\vol(G) \le (1 + \eps_0) \left( 1 - \frac{1}{\lambda^d} \right)\vol(\cW),
\end{equation}
\begin{equation} \label{eq:G_contains_small_neighborhoods}
\left(\cW \setminus \frac{1}{\lambda} \cW\right) + U_\delta \subseteq G.
\end{equation}

Denote $\ell_T = \# \left( TD_0 \cap \Lambda(G, \cL) \right)$. By \eqref{eq:irreducible_cut_and_project_density} when taking the window $G$, for large enough $T$, we have:
\begin{equation} \label{eq:initial_l_T_bound}
    \ell_T \le (1 + \eps_0) \cdot \frac{T^d \vol(G)}{\covol(\cL)} \ .
\end{equation}
Also denote the points in $TD_0 \cap \Lambda(G, \cL)$ by $(a_1^1, ...,a_d^1),...,(a_1^{\ell_T},...,a_d^{\ell_T})$.

Now order the prime ideals of $\cO_K$ according to their norms, $P_1,P_2,...$ such that $N(P_1) \le N(P_2) \le ...$. By the Chinese remainder theorem for prime ideals (see \cite[Appendix A]{marcus2018number}),
for every $1 \le i \le d$ we can find $x_i^0 \in \cO_K$, such that for all $1 \le j \le \ell_T$, we have $P_j | (x_i^0 + a_i^j)$. Equivalently, we have that:
\begin{equation} \label{eq:x0_has_large_invisible_neighborhood}
    P_j | (\bm{x}^0 + \bm{a}^j) \ \ \forall j \in \{1,...\ell_T\},
\end{equation}
with the same definition of dividing points as before.

Let $M_T > 0$ be a generator of the ideal $P_1 \cdot ... \cdot P_{\ell_T}$ (we can find a generator as $\cO_K$ is a principal ideal domain).

Denote the constant from \Cref{lem:lattice_point_in_rectangle} by $C > 0$. By this lemma, for each $1 \le i \le d$, we can find $(b_i, \sigma(b_i)) \in \cL_K$ for which:
\[
    \begin{dcases}
        &|M_T b_i + x_i^0| \le \frac{C N(M_T)}{\delta}
        ,\\
        &|\sigma(M) \sigma(b_i) + \sigma(x_i^0)| \le \delta,
    \end{dcases}
\]
as for each $1 \le i \le d$, this is a rectangle with volume $4C$, which is bigger than $C$. So denote $\bm{x} = (x_1^0 + M_T b_1,..., x_d^0 + M_T b_d)$, and we claim this is the $\bm{x}$ we want for \Cref{thm:large_holes_distance} (1). We have
\begin{equation} \label{eq:initial_x_norm_bound}
    \|\bm{x}\| \le \frac{C N(M_T)}{\delta},
\end{equation}
\begin{equation}\label{eq:x_has_small_sigma}
    \sigma(\bm{x}) \in U_\delta.
\end{equation}
with respect to the $L^\infty$ norm. By \eqref{eq:x0_has_large_invisible_neighborhood}, as $P_j | M_T$ for all $j\in\{1,...,\ell_T\}$, $\bm{x}$ also satisfies:
\begin{equation} \label{eq:x_has_large_invisible_neighborhood}
    P_j | (\bm{x} + \bm{a}^j) \ \ \forall j \in \{1,...\ell_T\}.
\end{equation}

\begin{lem} \label{lem:x_plus_TD_is_a_hole}
    For $\bm{x}$ defined above, we have:
    \[
    (\bm{x} +TD_0) \bigcap \Lambda_{vis}(\cW,\cL) = \es.
    \]
\end{lem}
\begin{proof}
    Assume by contradiction that $\bm{y} \in (\bm{x} + TD_0) \cap \Lambda_{vis}(\cW,\cL)$. By \Cref{prop:h_cond_visible_description}, we must have that no prime ideal divides $\bm{y}$, and $\sigma(\bm{y}) \in (\cW \setminus \frac{1}{\lambda} \cW)$. By \eqref{eq:x_has_small_sigma}, and \eqref{eq:G_contains_small_neighborhoods}, this also means that $\sigma(\bm{y} - \bm{x}) = \sigma(\bm{y}) - \sigma(\bm{x}) \in G$. As $\bm{y} \in \bm{x} + TD_0$, we get $\bm{y} - \bm{x} \in TD_0$. Together we have $\bm{y}-\bm{x} \in \left( TD_0 \cap \Lambda((A \setminus B), \cL) \right)$, thus there exists $1 \le j \le \ell_T$ with $\bm{y}-\bm{x} = \bm{a}^j$. So by \eqref{eq:x_has_large_invisible_neighborhood}, we get $P_j | \bm{x} + (\bm{y} - \bm{x}) = \bm{y}$, a contradiction to the visibility of $\bm{y}$.
\end{proof}

Now we bound the distance of $\bm{x}$ from the origin.
We have that $N(M_T) = \prod_{j=1}^{\ell_T}N(P_j)$. By the prime ideal theorem for number fields proved in \cite{Landau1903},
we have that $\frac{N(P_j)}{j \log{j}} \xrightarrow{j \rightarrow \infty} 1$. Thus we have that for large enough $T$:
\begin{equation}\label{eq:M_T_norm_bound}
    N(M_T)=\prod_{j=1}^{\ell_T}N(P_j) \le \exp((1 + \eps_0) \ell_T \log{\ell_T}).
\end{equation}

And also by \eqref{eq:initial_l_T_bound}:
\begin{equation} \label{eq:l_T_upper_bound}
    \ell_T \le (1 + \eps_0) \frac{T^d \vol(G)}{\covol(\cL)}
    \overset{\text{\eqref{eq:G_vol_bound}}}{\le} (1 + \eps_0)^2 \frac{\left(1 - \frac{1}{\lambda^d} \right) \vol(\cW)}{\covol(\cL)} T^d.
\end{equation}

Denote $R = \frac{\left(1 - \frac{1}{\lambda^d} \right) \vol(\cW)}{\covol(\cL)}$. We have that if $T$ is large enough, then $(1 + \eps_0)^2 R \le T^{\eps_0 d}$. Thus for large enough $T$, we have:
\begin{equation} \label{eq:log_l_T_bound}
    \log{\ell_T} \le (1 + \eps_0) d \log{T}.
\end{equation}

So combining the bounds, we get:
\[
    \begin{aligned}
        & \|x\| \overset{\text{\eqref{eq:initial_x_norm_bound}}}{\le} \frac{C N(M_T)}{\delta} \overset{\text{\eqref{eq:M_T_norm_bound}}}{\le} \exp{\left( (1 + \eps_0) \ell_T \log{\ell_T} + \log{C} + \log{\frac{1}{\delta}} \right)} \overset{\text{\eqref{eq:l_T_upper_bound},\eqref{eq:log_l_T_bound}}}{\le} \\
        & \le \exp{\left( (1 + \eps_0) \cdot (1 + \eps_0)^2 R T^d \cdot(1+ \eps_0) d \log{T} + \log{C} + \log{\frac{1}{\delta}} \right)} \le \\
        &
        \overset{\textit{For T large enough}}{\le} \exp{\left( (1 + \eps_0)(1 + \eps_0)^4 R T^d d \log{T} \right)} = \\
        & =\exp{\left( (1 + \eps) \frac{\left(1 - \frac{1}{\lambda^d} \right) \vol(\cW)}{\covol(\cL)} d T^d \log{T} \right)}.
    \end{aligned}
\]
Together with \Cref{lem:x_plus_TD_is_a_hole}, this proves the claim.
\end{proof}

\subsection{Lower bound on the distance,  \Cref{thm:large_holes_distance} (2)}
\begin{lem} \label{lem:point_counting_L_g_for_grid}
    Let $D$ be a nice averaging set, $\cW$ a nice window, and $r > 0$. Then there exists a constant $C = C(D, \cW, r)$ such that for any $0 \ne I \vartriangleleft \cO_K$, $T \ge N(I)$, and $\bm{y}_1, \bm{y}_2 \in \bR^d$, we have:
    \[
    \left| \# \left(\Lambda(r \cW + \bm{y}_2,\cL_I) \cap (\bm{y}_1 + TD)\right) - \frac{\vol(TD) \cdot \vol(\cW)}{\covol(\cL)} \cdot \frac{r^d}{N(I)^d} \right| \le C \left( \frac{T}{N(I)} \right)^{d-\frac{1}{2}}.
    \]
\end{lem}
\begin{proof}
    This is again a generalization of \cite[Lemma 5.4]{gringlaz2025errorboundsvisiblepoints}. Here we also add that we can change the location of the window and the measuring set, and it won't affect the bounds. We also claim the bound is for the lattice $\cL_I$ instead of $a_\fkg \cL$. The exact same proof still works, as the diameters and volumes calculated in the proof stay the same and not affected by $\bm{y}_1,\bm{y}_2$, and \cite[Proposition 5.2]{gringlaz2025errorboundsvisiblepoints} (which is from \cite{Schmidt1995}) doesn't depend on moving the convex set in the space. Also working with $\cL_I$ instead of $a_\fkg \cL$ doesn't matter, as we assume $K$ has H-condition, and in particular $\cO_K$ is a principle ideal domain, so for every ideal $I$ we have $\fkg \in \cO_K$ for which $I = (\fkg)$, thus $\cL_I = a_\fkg \cL$, and $N(\fkg) = N(I)$.
\end{proof}
\begin{proof} [Proof of \Cref{thm:large_holes_distance} (2)]
As $\cW$ is bounded, take $A > 1$ constant such that $\cW \subseteq [-A,A]^d$. As before, order the prime ideals of $\cO_K$ according to their norms, $P_1,P_2,...$ such that $N(P_1) \le N(P_2) \le ...\ $. For convenience, denote $\cW_\lambda := \cW \setminus \frac{1}{\lambda} \cW$.
    Denote $k,\ell \in \bN$ to be the maximal indices for which
    \begin{equation} \label{eq:k_and_l_definition}
    N\left( \prod_{i=1}^k{P_i} \right) \le \sqrt{T},\ \ \ N(P_\ell) \le 2AT,
    \end{equation}
    both hold. As we know from \Cref{prop:h_cond_visible_description}, the set of visible points is:
    \[
    \Lambda_{vis}(\cW,\cL) = \left\{ (x_1,..., x_d) \in \Lambda_{\star}\left(\cW_\lambda, \cL\right) :  \, \gcd(x_1,..., x_d)=(1) \right\}.
    \]
    We are given that $(\bm{x} +TD_0) \bigcap \Lambda_{vis}(\cW,\cL) = \es$. Thus each point in $\Lambda_{\star}\left(\cW_\lambda, \cL\right) \cap (\bm{x} + TD_0)$ must be divisible by some prime ideal. The main idea will be to bound from below how many points in $\Lambda_{\star}\left(\cW_\lambda, \cL\right) \cap (\bm{x} + TD_0)$ are not divisible by any of $P_1,...,P_k$ all together using the Chinese remainder theorem in \Cref{lem:bound_small_ideals_count_crt}, then to bound from above how many of those points are divisible by each of $P_{k+1},...,P_\ell$ individually in \Cref{lem:middle_ideals_union_bound}, in order to get that at least $\eps T^d$ (for $\eps > 0$ a constant) points in $\Lambda_{\star}\left(\cW_\lambda, \cL\right) \cap (\bm{x} + TD_0)$ must be divisible by $P$ with $N(P) > 2AT$, and finally getting a lower bound on $\|\bm{x}\|$ from those primes of large norm.

    For a nice averaging set $D$, use \Cref{lem:point_counting_L_g_for_grid} for $r = 1, \frac{1}{\lambda}$ and our window $\cW$, and denote $C_D$ the sum of the constants for the two cases. So for any $I,T,\bm{y}_1, \bm{y}_2$ as in the statement (with the main requirement $N(I) \le T$) we have:
    \begin{equation} \label{eq:W_lambda_point_counting}
        \begin{aligned}
            &\left| \# \left(\Lambda\left(\cW_\lambda + \bm{y}_2,\cL_I\right) \cap (\bm{y}_1 + TD)\right) - \frac{\vol(TD)\vol(\cW_\lambda)}{\covol(\cL) N(I)^d} \right| \le \\
            &\le C_D \left( \frac{T}{N(I)} \right)^{d-\frac{1}{2}},
        \end{aligned}
    \end{equation}
    by taking the difference between the lemma for $r = 1, \frac{1}{\lambda}$. We will use \Cref{eq:W_lambda_point_counting} with $D = D_0$ in \Cref{lem:bound_small_ideals_count_crt}, and with $D = (2A)^{-1}D_0$ in \Cref{lem:middle_ideals_union_bound}.
\begin{lem} \label{lem:bound_small_ideals_count_crt}
    There exists a constant $C_1 > 0$, such that for large enough $T > 1$ and $k \in \bN$ from \eqref{eq:k_and_l_definition}, we have:
    \[
    \# \left(\Lambda\left(\cW_\lambda,\cL \setminus \left( \bigcup_{i=1}^{k}{\cL_{P_i}}\right)\right) \cap (\bm{x} + TD_0) \right) \ge T^d \left(\theta(\Lambda_{vis}(\cW, \cL)) - C_1 \cdot T^{-\frac{1}{4}}\right).
    \]
\end{lem}
\begin{proof}
    Denote:
    \[
    I = \prod_{i=1}^k{P_i},\quad s = \prod_{i=1}^k(N(P_i)^d-1).
    \]
    Using the Chinese remainder theorem for prime ideals, take $\bm{y}^1,\bm{y}^2,...,\bm{y}^s \in \cO_K^d$ to be representatives of the points in $\left( \cO_K / I \right) ^d$, which are not divisible by any of $P_1,...,P_k$. In other words, for any $\bm{y} \in \cO_K^d$, there exists $j \in \{1,...,s\}$ such that $I | (\bm{y} - \bm{y}^j)$ if and only if, for all $i = 1,...,k$, we have that $P_i \not| \ \bm{y}$. And also if such $j$ exists, it is unique. Denote also $\bm{z}^j = (\bm{y}^j , \sigma(\bm{y}^j))$. So we can write the disjoint union:
    \[
    \cL \setminus \left( \bigcup_{i=1}^{k}{\cL_{P_i}}\right) = \biguplus_{j=1}^{s}{\left(\cL_I + \bm{z}^j\right)}.
    \]
    Note that we required $N(I) = N(\prod_{i=1}^k{P_i}) \le \sqrt{T} \le T$, so \eqref{eq:W_lambda_point_counting} requirement is satisfied, so using it with $D = D_0$, we have:
    \[
    \begin{aligned}
        &\# \left(\Lambda\left(\cW_\lambda,\cL \setminus \left( \bigcup_{i=1}^{k}{\cL_{P_i}}\right)\right) \cap (\bm{x} + TD_0) \right) = \\
        & = \sum_{j=1}^s \# \left(\Lambda\left(\cW_\lambda,\cL_I + \bm{z}^j)\right) \cap (\bm{x} + TD_0) \right) = \\
        & = \sum_{j=1}^s \# \left(\Lambda\left(\cW_\lambda - \sigma(\bm{y}^j),\cL_I\right) \cap (\bm{x} - \bm{y}^j + TD_0) \right) \overset{\eqref{eq:W_lambda_point_counting}}{\ge} \\
        & \ge \sum_{j=1}^s \left(\frac{\vol(TD_0)\vol(\cW_\lambda)}{\covol(\cL) N(I)^d} - C_{D_0} \left( \frac{T}{N(I)} \right)^{d-\frac{1}{2}} \right) = \\
        & = \frac{s}{N(I)^d} \left(\frac{\vol(TD_0)\vol(\cW_\lambda)}{\covol(\cL)} - C_{D_0} T^{d-\frac{1}{2}}\left(N(I)\right)^{\frac{1}{2}} \right).
    \end{aligned}
    \]
    Now noting that $s = \prod_{i=1}^k(N(P_i)^d-1)$, we have:
    \[
    \frac{s}{N(I)^d} = \prod_{i=1}^{k}{\left(1 - \frac{1}{N(P_i)^d}\right)} \ge \prod_{i=1}^{\infty}{\left(1 - \frac{1}{N(P_i)^d}\right)} = \frac{1}{\zeta_K(d)}.
    \]
    So also bounding $N(I) \le \sqrt{T}$, and substituting $\theta(\Lambda_{vis}(\cW, \cL)) = \frac{\vol(\cW_\lambda)}{\zeta_K(d)\covol(\cL)}$ (as proved in \cite{hammarhjelm2019densityminimalgapvisible}, see \Cref{rmk:H_condition_constants}) in total we have at least:
    \[
    T^d \left(\theta(\Lambda_{vis}(\cW, \cL)) - \frac{C_{D_0}}{\zeta_K(d)} \cdot T^{-\frac{1}{4}}\right)
    \]
    points in $\Lambda_{\star}\left(\cW_\lambda, \cL\right) \cap (\bm{x} + TD_0)$ which are not divisible by any of $P_1,...,P_k$. So taking $C_1 = \frac{C_{D_0}}{\zeta_K(d)}$ gives us what we wanted.
\end{proof}
\begin{lem} \label{lem:middle_ideals_union_bound}
    For any $\eps > 0$, for large enough $T$ and $k,\ell$ from \eqref{eq:k_and_l_definition}, we have:
    \[
    \sum_{i = k+1}^{\ell}{\# \left(\Lambda\left(\cW_\lambda,\cL_{P_i} \right) \cap (\bm{x} + TD_0) \right)} \le \eps T^d.
    \]
\end{lem}
\begin{proof}
    We use \eqref{eq:W_lambda_point_counting} but with $D = (2A)^{-1}D_0$ instead of $D_0$. Define $C_2 = (2A)^{d - \frac{1}{2}}C_D$ and denote $T' = 2AT$, which we will use in \eqref{eq:W_lambda_point_counting} instead of $T$. Now for each $k < i \le \ell$ taking $I = P_i$, we have $N(I) \le N(P_\ell) \le 2AT = T'$, and take $\bm{y_1} = \bm{x}, \bm{y}_2 = 0$, so we get:
    \[
    \# \left(\Lambda\left(\cW_\lambda,\cL_{P_i}\right) \cap (\bm{x} + TD_0)\right) \le \frac{\vol(TD_0)\vol(\cW_\lambda)}{\covol(\cL) N(P_i)^d} + C_2 \left( \frac{T}{N(P_i)} \right)^{d-\frac{1}{2}}.
    \]
    So summing over $k < i \le \ell$, we have:
    \[
    \begin{aligned}
        &\sum_{i=k+1}^{\ell}{\left(\frac{\vol(TD_0)\vol(\cW_\lambda)}{\covol(\cL) N(P_i)^d} + C_2 \left( \frac{T}{N(P_i)} \right)^{d-\frac{1}{2}}\right)} \le \\
        & \le T^d \frac{\vol(\cW)}{\covol(\cL)} \left(\sum_{i=k+1}^{\infty}{\frac{1}{N(P_i)^d}} + \frac{C_3}{\sqrt{T}}\sum_{i=k+1}^{\infty}{\frac{1}{N(P_i)^{d-\frac{1}{2}}}}\right),
    \end{aligned}
    \]
    where $C_3 = C_2 \frac{\covol(\cL)}{\vol(\cW)}$. So as $T \rightarrow \infty$, also $k \rightarrow \infty$, and $\sum_{i=k+1}^{\infty}{\frac{1}{N(P_i)^d}}$ is the tail of a converging series (the entire series is bounded by $\sum_{0 \ne I \vartriangleleft \cO_K}{\frac{1}{N(I)^d}} = \zeta_K(d)$), so it goes to zero. The second term in the brackets goes to zero for similar reasons. So for any $\eps > 0 $, for $T$ large enough, at most $\eps T^d$ are blocked by $P_i$ for $k < i \le \ell$.
\end{proof}
    From here throughout the rest of the proof, denote $\eps = \frac{1}{4}\theta(\Lambda_{vis}(\cW, \cL))$, and we can take $T$ large enough so \Cref{lem:bound_small_ideals_count_crt,lem:middle_ideals_union_bound} hold for it, as $\eps$ is constant.
    For $T$ large enough, note we also have $C_1 \cdot T^{-\frac{1}{4}} \le \eps$, so now we have that:
    \[
    \begin{aligned}
        &\# \left(\Lambda\left(\cW_\lambda,\cL \setminus \left( \bigcup_{i=1}^{\ell}{\cL_{P_i}}\right)\right) \cap (\bm{x} + TD_0) \right) \ge \\
        & \# \left(\Lambda\left(\cW_\lambda,\cL \setminus \left( \bigcup_{i=1}^{k}{\cL_{P_i}}\right)\right) \cap (\bm{x} + TD_0) \right) - \sum_{i = k+1}^{\ell}{\# \left(\Lambda\left(\cW_\lambda,\cL_{P_i} \right) \cap (\bm{x} + TD_0) \right)} \ge \\
        & \ge T^d \left(\theta(\Lambda_{vis}(\cW, \cL)) - \eps\right) - \eps T^d = 2\eps T^d.
    \end{aligned}
    \]
    So at least $2\eps T^d$ points in $\Lambda\left(\cW_\lambda, \cL\right) \cap (\bm{x} + TD_0)$ are divisible only by primes $P$ with $N(P) > 2AT$.

    \begin{lem} \label{lem:large_norm_primes_divide_one_point}
    For any prime ideal $0 \ne P \vartriangleleft \cO_K$ with $N(P) > 2AT$, it divides at most one point in $\Lambda\left(\cW, \cL\right) \cap (\bm{x} + TD_0)$.
    \end{lem}
    \begin{proof}
        Assume otherwise, so we have $\bm{x}^1, \bm{x}^2 \in \Lambda\left(\cW, \cL\right) \cap (\bm{x} + TD_0)$, with $P | \bm{x}^1, \bm{x}^2$. So $\bm{x}^1,\bm{x}^2 \in \bm{x} +[0,T]^d$, and $\sigma(\bm{x}^1), \sigma(\bm{x}^2) \in \cW \subseteq [-A,A]^d$. Thus $P | (\bm{x}^1 - \bm{x}^2)$, but also, $\bm{x}^1-\bm{x}^2 \in [-T,T]^d$, and $\sigma(\bm{x}^1 - \bm{x}^2) \in [-2A,2A]^d$. Looking at the $i$-th coordinate, we have that $|x_i^1-x_i^2| \le T$, and $|\sigma(x_i^1 - x_i^2)| \le 2A$, so $|N(x_i^1-x_i^2)| \le 2AT < N(P)$. But $P | (x_i^1-x_i^2)$, so $N(P) | N(x_i^1-x_i^2)$, thus we must have that $|N(x_i^1 - x_i^2)| = 0 $ and $x_i^1 = x_i^2$. This is for every coordinate, so $\bm{x}^1 = \bm{x}^2$.
    \end{proof}

    Now denote $\pi_i : \bR^d \rightarrow \bR$ for $i \in \{1,...,d\}$ to be the projection to the $i$'th coordinate. Define the finite set:
    \[
    F = \pi_1\left(\Lambda(\cW,\cL) \cap (\bm{x} +TD_0)\right) \setminus \{0\} \subseteq \cO_K.
    \]
    And let us look at the product:
    \[
    R = \left(\prod_{y \in F}{y}\right) \vartriangleleft \cO_K.
    \]
    Note now that $R$ is divisible by all prime ideals that divide points in $\Lambda(\cW,\cL) \cap (\bm{x} +TD_0)$ with non-zero first coordinate. We would like to use this fact to bound $N(R)$ from below, so first let us bound from above the number of points with zero first coordinate:
    
    As $\cL_K = \{ (x, \sigma(x)) : x \in \cO_K \}$ is a lattice, there exists $C_4 >0 $ depending only on $K$ such that for any $A,T \ge 1$ and $t \in \bR$ we have:
    \begin{equation} \label{eq:bound_rectangle_L_K_point_count}
        \#\left(\cL_K \cap \left( (t + [0,T]) \times [-A,A] \right) \right) \le C_4 AT.
    \end{equation}
    So we have:
    \[
    \begin{aligned}
        &\#\left(\pi_1^{-1}(0) \cap \Lambda(\cW,\cL) \cap (\bm{x} +TD_0) \right) \le \\
        &\le \#\left(\pi_1^{-1}(0) \cap \Lambda([-A,A]^d,\cL) \cap (\bm{x} +TD_0) \right) \le (C_4 AT)^{d-1}.
    \end{aligned}
    \]
    So for large enough $T > 0$, there are at most $\eps T^d$ such points. So now we are left with at least $2\eps T^d - \eps T^d = \eps T^d$ points in $\Lambda(\cW,\cL) \cap (\bm{x} +TD_0)$ which are divisible by prime ideals of norm larger than $2AT$, and have non-zero first coordinate. Each prime dividing one of them divides $R$ as well. By \Cref{lem:large_norm_primes_divide_one_point}, each of those prime ideals divides at most one point in $\Lambda(\cW,\cL) \cap (\bm{x} +TD_0)$, thus there are at least $\eps T^d$ different prime ideals with norm larger than $2AT$ dividing $R$. Also $R \ne (0)$, so we can bound $N(R)$ from below by the product of the norms of the primes:
    \begin{equation} \label{eq:projection_norm_lower_bound}
        N(R) \ge (2AT)^{\eps T^d}.
    \end{equation}
    
    On the other hand, by \eqref{eq:bound_rectangle_L_K_point_count}, we also have:
    \[
    \begin{aligned}
        &\#F \le \#\left(\pi_1\left(\Lambda\left([-A,A]^d, \cL\right) \cap (\bm{x} + TD_0)\right)\right) \le \\
        & \#\left(\cL_K \cap \left( (\pi_1(\bm{x}) + [0,T]) \times [-A,A] \right) \right) \le C_4 AT.
    \end{aligned}
    \]
    Note for every $y \in F$, we have $|y| \le |\pi_1(\bm{x})| + T$ and $|\sigma(y)| \le A$. Thus we can bound $N(R)$ from above:
    \begin{equation} \label{eq:projection_norm_upper_bound}
        \begin{aligned}
            &N(R) = N\left(\prod_{y \in F}{y}\right) = \prod_{y \in F}{|y\sigma(y)|} \le \\
            & \le (A(|\pi_1(\bm{x})| + T))^{\#F} \le (A(|\pi_1(\bm{x})| + T))^{C_4 AT}.
        \end{aligned}
    \end{equation}
    
    Thus combining \cref{eq:projection_norm_lower_bound,eq:projection_norm_upper_bound}, we have:
    \[
    (A(|\pi_1(\bm{x})| + T))^{C_4 AT} \ge (2AT)^{\eps T^d}.
    \]
    Moving terms around, we get:
    \[
    |\pi_1(\bm{x})| \ge \frac{1}{A}(2AT)^{\frac{\eps}{C_4A} T^{d-1}} - T = \exp{\left( \frac{\eps}{C_4A} T^{d-1} (\ln(T)+ \ln(2A)) - \ln(A)\right)} - T.
    \]
    Taking some $0< c < \frac{\eps}{C_4A}$ constant, for $T$ large enough, this is bigger than $\exp{\left( c T^{d-1} \log(T)\right)}$, thus $|\pi_1(\bm{x})| \ge \exp{\left( c T^{d-1} \log(T)\right)}$, and:
    \[
    \|\bm{x}\| \ge \exp{\left( c T^{d-1} \log(T)\right)},
    \]
    as we wanted.
    
\end{proof}

\section*{Acknowledgments}

This paper is a part of the author’s M.Sc. thesis at Tel Aviv University under the supervision of Prof. Barak Weiss. The support of ISF grant 2021/24 is gratefully acknowledged. The author thanks Barak Weiss for helpful discussions, directions, and for comments on drafts of the paper which vastly helped improving it.
The author also thanks Omri Peer for helpful discussions.

\subsection*{AI disclosure statement}
AI was used only for some of the LaTeX throughout the paper, the wording of some of the proofs, and for finding several of the appropriate references. It wasn't used for the research and proofs themselves.

\bibliographystyle{plain}
\bibliography{references}

\begin{thebibliography}{10}

\bibitem{BaakeGrimm2013}
Michael Baake and Uwe Grimm.
\newblock {\em Aperiodic Order. Volume 1: A Mathematical Invitation}, volume 149 of {\em Encyclopedia of Mathematics and its Applications}.
\newblock Cambridge University Press, Cambridge, 2013.

\bibitem{baake2000diffraction}
Michael Baake, Robert~V Moody, and Peter~A.B. Pleasants.
\newblock Diffraction from visible lattice points and {$k$}th power free integers.
\newblock {\em Discrete Mathematics}, 221(1-3):3--42, 2000.

\bibitem{BeenkerNonPeriodicTilings1982}
F.P.M. Beenker.
\newblock {\em Algebraic theory of non-periodic tilings of the plane by two simple building blocks : a square and a rhombus}.
\newblock EUT report. WSK, Dept. of Mathematics and Computing Science. Eindhoven University of Technology, 1982.

\bibitem{DEBRUIJN198153}
N.G. {de Bruijn}.
\newblock Algebraic theory of penrose's non-periodic tilings of the plane. ii.
\newblock {\em Indagationes Mathematicae (Proceedings)}, 84(1):53--66, 1981.

\bibitem{Erdos_1958}
Paul Erdős.
\newblock On an elementary problem in number theory.
\newblock {\em Canadian Mathematical Bulletin}, 1(1):5–8, 1958.

\bibitem{gringlaz2025errorboundsvisiblepoints}
Ilya Gringlaz, Rishi Kumar, and Barak Weiss.
\newblock On the error bounds for visible points in some cut-and-project sets, 2025.
\newblock arXiv:2505.02631 [math.DS].

\bibitem{hammarhjelm2019densityminimalgapvisible}
Gustav Hammarhjelm.
\newblock The density and minimal gap of visible points in some planar quasicrystals.
\newblock {\em Discrete Mathematics}, 345(12):113074, 2022.

\bibitem{hof1998uniform}
Anton Hof.
\newblock Uniform distribution and the projection method.
\newblock In Ji\v{r}\'{i} Patera, editor, {\em Quasicrystals and Discrete Geometry}, volume~10 of {\em Fields Institute Monographs}, pages 201--206. American Mathematical Society, Providence, RI, 1998.

\bibitem{kumar2026simultaneous}
Rishi Kumar and Carlos Ospina.
\newblock Simultaneous visibility in cut-and-project sets.
\newblock {\em arXiv preprint arXiv:2610.28119}, 2026.

\bibitem{Landau1903}
Edmund Landau.
\newblock Neuer beweis des primzahlsatzes und beweis des primidealsatzes.
\newblock {\em Mathematische Annalen}, 56(4):645--670, 1903.

\bibitem{marcus2018number}
D.A. Marcus.
\newblock {\em Number Fields}.
\newblock Universitext. Springer-Verlag New York, 2018.

\bibitem{marklof2014visibilitydirectionsquasicrystals}
Jens Marklof and Andreas Str{\"o}mbergsson.
\newblock Visibility and directions in quasicrystals.
\newblock {\em International Mathematics Research Notices}, 2015(15):6588--6617, 2015.

\bibitem{Marklof_2014}
Jens Marklof and Andreas Strömbergsson.
\newblock Free path lengths in quasicrystals.
\newblock {\em Communications in Mathematical Physics}, 330(2):723–755, March 2014.

\bibitem{Meyer1972}
Yves Meyer.
\newblock {\em Algebraic Numbers and Harmonic Analysis}, volume~2 of {\em North-Holland Mathematical Library}.
\newblock North-Holland Publishing Company, Amsterdam, 1972.

\bibitem{Moody_2002}
Robert~V. Moody.
\newblock Uniform distribution in model sets.
\newblock {\em Canadian Mathematical Bulletin}, 45(1):123–130, 2002.

\bibitem{rudnick2026classification}
Zeev Rudnick.
\newblock The classification of real quadratic fields which satisfy hammarhjelm's condition.
\newblock {\em arXiv preprint arXiv:2606.31544}, 2026.

\bibitem{schlottmann1998cut}
Martin Schlottmann.
\newblock Cut-and-project sets in locally compact abelian groups.
\newblock {\em Quasicrystals and discrete geometry}, 10:247--264, 1998.

\bibitem{Schmidt1995}
Wolfgang~M. Schmidt.
\newblock Northcott's theorem on heights {II}. {T}he quadratic case.
\newblock {\em Acta Arithmetica}, 70(4):343--375, 1995.

\bibitem{simmons2016danzer}
David Simmons and Yaar Solomon.
\newblock A danzer set for axis parallel boxes.
\newblock {\em Proceedings of the american mathematical society}, 144(6):2725--2729, 2016.

\end{thebibliography}

\end{document}